%% file: quad-n-quat-algs.tex
\documentclass[11pt,letterpaper]{article}

\pdfoutput=1

\usepackage{hyperref}

\input{extended-preamble.tex}

\addtotextwidth{2cm}
\addtotextheight{2.5cm}

\input{preamble.tex}

\begin{document}

\title{\vspace{-6.5ex}
	\LARGE\bfseries Quadratic extension algebras\\ and quaternion algebras\\
	over fields (of characteristic not $\boldsymbol{2}$)\\[2ex]}

\author{\large\bfseries France Dacar, Jo\v{z}ef Stefan Institute\\
	{\large\ttfamily France.Dacar@ijs.si}}

\date{\large\bfseries 
	March 28, 2013\\
	Last edited August 7, 2016\\[10ex]}

\maketitle

\begin{abstract}
The paper presents a classification of quadratic extension algebras,
	also known as algebras of degree~2,
as well as several characterizations of quaternion algebras over a~field
	(of characteristic not~2).
The presentation is not restricted to finite-dimensional algebras.
The `pure calculus' on a quaternion algebra is introduced;
	it~generalizes the `vector calculus' of the Hamilton's quaternions,
and is instrumental in establishing, in a~clear and simple way,
the isomorphism between
	the group consisting of all automorphisms and all anti-automorphisms of a~quaternion algebra
	and the orthogonal group of the norm form on the subspace of pure quaternions.
The~paper concludes with a~glimpse of quaternion algebras over an integral domain
	(of characteristic not~2),
and of their facility in studying ternary quadratic forms over an integral domain.
\end{abstract}

\hypersetup{
  pdftitle={Quadratic extension algebras and quaternion algebras over fields (of characteristic not 2)},
  pdfauthor={France Dacar},
  pdfkeywords={quadratic extension algebra, algebra of degree 2, Clifford algebra, quaternion algebra}
}

\section{Prologue}
\label{sec:prologue}

\inskip

Tsit Yuen Lam \cite{LamIQFF},
	the end of section III.5, ``$\tinysp$Characterizations of Quaternion Algebras'':
\begin{quotation}
The second characterization of quaternion algebras
is in terms of the canonical involution ``$\tinysp$bar'' on such algebras.
Note that this is an involution ``of the first kind'':
that is, it is an involution
that restricts to the identity map on the center $F$ of the algebra.

\txtskip

\noindent{\bf Theorem 5.2.}\: \textit{Let\/ $B\neq F$ be
a~finite-dimensional simple\/ $F\negdtinysp$-algebra with center\/ $F$
equipped with an\/~$F\negdtinysp$-algebra involution of the first kind\/ $x\mapsto\overbar{x}$
such~that\/ $x+\overbar{x}\in F$ and\/ $x\tinysp\overbar{x}\in F$ {\rm(}for all\/ $x\narrt\in B${\rm)}.
Then\/ $B$ is isomorphic to~a~quaternion algebra over\/~$F\negdtinysp$.}

\txtskip

Although this is a very nice characterization of quaternion algebras,
it will not be needed in the rest of this book.
Therefore, we'll leave its proof as an exercise to the reader.
The main idea of the proof is that,
given the~properties of the involution $x\mapsto\overbar{x}$,
every element $x\in B$ satisfies a quadratic equation over the center~$F$;
namely
\begin{equation*}
x^2 -x\tinysp(x\narrt+\overbar{x})+x\tinysp\overbar{x}=0
	\qquad (x\narrt+\overbar{x},\,x\tinysp\overbar{x}\in F).
\end{equation*}
\end{quotation}
From all properties of the $F\negdtinysp$-algebra $B$ mentioned in the quoted text
we shall single out the property
``every element $x\in B$ satisfies a quadratic equation over~$F\tinysp$'',
and shall try to find out how far it gets us.
The idea is to gather enough information about such algebras,
so that by imposing on them some additional conditions,
preferably as weak as possible, 
we will be able to characterize quaternion algebras.

\txtskip

We kicked this particular stone that rested peacefully at the top of the hill,
and it started rolling and tumbling 
down the slope.
We are going to run after it,
because we are curious to see where 
it will come to rest.

Ready?
Here we go.

\pagebreak[3]\txtskip

\section{Some conventions, terminology, and notation}
\label{sec:prelims}

\medskip

\emph{All fields mentioned in this text
are assumed to have the characteristic different from~two.
In particular, $K\negtinysp$ will always denote a~field with\/ $\fieldchar{K}\negtinysp\neq 2$. %
Unless specified otherwise, every algebra is over a~field
and is silently assumed to be associative, unital,
and nontrivial (\ie, it does not consist of the zero alone, which also serves as~the~unity).}

\txtskip

We will encounter some algebras that do not conform to the default assumptions.
The \notion{zero algebra} (over some field) is the set $\set{0}$ equipped with the evident operations;
the solitary element~$0$ of the zero algebra is the neutral element both for addition and multiplication.
A~\notion{zero$\tinysp$-product algebra} is a vector space $A$ equipped with the all-zero multiplication,
that is, $a\tinysp b=0$ for all $a,\tinysp b\in A\tinysp$;
a~zero$\tinysp$-product algebra is associative,
and it is unital \iff\ it is the zero algebra.
A~\notion{zero$\tinysp$-square algebra} is an associative algebra~$A$
	in which $a^2=0$ for every $a\in A$;
it is unital \iff\ it is the zero algebra.
Two other (equivalent) descriptions of a~zero$\tinysp$-square algebra:
an~anti\-commutative associative algebra; an associative Lie algebra.

\txtskip

Let $A$ be an associative algebra, not necessarily with unity, over a~field $K\negtinysp$.
We~\notion{adjoin an identity element (unity)} to $A$
to obtain the algebra $K\negtinysp\oplus A\tinysp$;
writing $\pair{\lambda,a}\in K\negtinysp\oplus A$
as a~formal sum $\lambda\oplus a$,
the multiplication on $K\negtinysp\oplus A$ is defined by
$$(\lambda\narrt\oplus a)(\mu\narrt\oplus b)
	\Defeq \lambda\mu\oplus(\lambda b\narrt+\mu\tinysp a\narrt+ ab)~.$$
The algebra $K\negtinysp\oplus A$ is associative with the unity~$1\oplus 0$.%
\footnote{Adjoining an~identity element is a~universal construction.
Let $K\negtinysp$ be a field,
let $\cat{Alg}_K$ be the category of associative unital $K\negdtinysp$-algebras
with morphisms the functions preserving addition, multiplication by scalars, and multiplication,
and sending the unity to the unity;
also, let $\cat{Alg}'_K$ be the category of associative $K\negdtinysp$-algebras
whose morphisms are functions preserving addition, multiplication by scalars, and multiplication.
Let~$\cat{Alg}_K\to\nolinebreak\cat{Alg}'_K$ be the functor which `forgets the unity';
the left adjoint $A\mapsto K\negtinysp\oplus A$ to this functor
does the job of adjoining identity elements,
and for each algebra $A$ in $\cat{Alg}'_K$
the mapping $A\narrt\to K\negtinysp\narrt\oplus A \narrt: a\narrt\mapsto 0\oplus a$
is~a~universal arrow in $\cat{Alg}'_K$ from~$A$ to the forgetful functor.}
If the algebra $A$ already possesses a~unity~$1_{\negtinysp A}$,
then $0\oplus 1_{\negtinysp A}$ is an idempotent of the unital algebra $K\negtinysp\oplus A$
that is different from its unity.

\pagebreak[3]
\txtskip

Whenever $S$ is (the underlying set of) a~structure
containing a~distinguished element called ``$\tinysp$zero''
and denoted by~$0$,
and $X$ is any subset of~$S$,
then we write $\punct{X}\defeq X\narrdt\setdiff\set{0}$.

If $S$ is a~structure,
part of which is a~multiplicatively written associative binary operation with a~neutral element
(and there is only one such operation in the structure),
which therefore endows $S$ with the structure of a~monoid,
then we denote this monoid by $S^{\tinysp\cdot}$ (note the small dot)
and by $S^\times$ the multiplicative group of invertible elements of~$S^{\tinysp\cdot}$.
In particular, if $S$ is a~ring, or an algebra,
then $S^\times$ denotes the group of units (multiplicatively invertible elements) of~$S$.

If $S$ is a~structure, part of which is an additively written binary operation
which makes~$S$ an abelian group
(examples of such structures are rings, algebras, modules,~\ldots),
then~we~denote this abelian group by $S^+\negdtinysp$.

Let $S$ be a~structure, part of which is a~multiplicatively written binary operation.
For every subset $X$ of~$S$ we shall denote by $\sqr{X}$
the set of squares $xx$ of all elements~$x$~of~$X$.
If~$S$ is a~structure for which we have also defined $\punct{X}$ for $X\subseteq S$,
or~$S^\times\negdtinysp$, or both,
then we agree that $\sqr\punct{X}$ means $\sqr(\punct{X})$,
and $\sqr S^\times$ means $\sqr(S^\times)$;
that is, $\punct{(\anon)}$ and $(\anon)^\times$ have precedence over~$\sqr(\anon)$.

\pagebreak[3]\bigskip

\section{Quadratic extension algebras}
\label{sec:quad-ext-algs}

\medskip

We define a~\notion{quadratic extension\/ $K\negdtinysp$-algebra}%
\footnote{Also known as a $K\negdtinysp$-algebra of degree~$2$.}
%
%
%
as a~non-zero associative $K\negdtinysp$-algebra $A$ with unity $1=1_A$,
where every $x\in A$ satisfies a~quadratic equation
$x^2+\alpha\tinysp x+\beta=0$ for some $\alpha,\tinysp\beta\in K\negtinysp$.
Let $A$ be a quadratic extension $K\negdtinysp$-algebra.
As~usual we identify $K1_A$ with~$K\negtinysp$ (we can do this because $1_A\neq 0$)
so that $K\negtinysp\subseteq A$,
and refer to the elements of $K\negtinysp$ as~\notion{scalars}.
We~set $\pure{A}\defeq\set{0}\union\set{u\in A\narrdt\setdiff K\negtinysp\suchthat u^2\in K\negtinysp}$,
call $\pure{A}$ the \notion{pure part} of~$A$,
and say that elements of $\pure{A}$ are the \notion{pure elements} of~$A$.

\thmskip

\begin{lemma}\label{lem:scalar-pure-split}
Every element\/ $x$ of a~quadratic extension\/ $K\negdtinysp$-algebra\/ $A$
has a~unique representation
of the form\/ $x=\alpha+u$ with\/ $\alpha\in K\negtinysp$ and\/ $u\in\pure{A}$.
\end{lemma}

\interskip

\begin{proof}
If $x=\alpha\in K\negtinysp$, then $x=\alpha+0$ is a~required representation,
which is unique: 
if~$x=\beta+u$ with $\beta\in K$ and $u\in\pure{A}$, then $u=\alpha-\beta\in K\inters\pure{A}=\set{0}$.

Suppose $x\notin K\negtinysp$.
Then $x^2+\beta\tinysp x+\gamma=0$ for some $\beta,\tinysp\gamma\in K\negtinysp$,
thus $\bigl(x+\thalf\beta\bigr)^{\negtinysp 2}=\tfrac{1}{4}\beta^{\tinysp2}-\gamma$,
which means that $x=\alpha+u$
with $\alpha=-\thalf\beta\in K\negtinysp$ and $u^2=\tfrac{1}{4}\beta^{\tinysp2}-\gamma\in K\negtinysp$;
since~$u\in K$ implies $x\in K$,
we must have $u\notin K\negtinysp$ because $x\notin K$, and hence $u\in\pure{A}$.
Uniqueness.
Suppose that also $x=\alphapr+\upr$ with $\alphapr\in K\negtinysp$ and $\upr\in\pure{A}$.
Then $\upr=(\alpha\narrt-\alphapr)+u$,
where $(\upr)^2=(\alpha\narrt-\alphapr)^2+u^2+2\tinysp(\alpha\narrt-\alphapr)\tinysp u\in K$
and $u^2\in K$,
thus $\delta\defeq(\alpha\narrt-\alphapr)\tinysp u\in K\negtinysp$,
so~we must have $\alpha-\alphapr=0$,
since otherwise $u = \delta/(\alpha\narrt-\alphapr)\in K$, contradicting $u\notin K\negtinysp$.
\end{proof}

\pagebreak[3]
\thmskip

We will occasionally introduce (or represent)
an element of a~quadratic extension $K\negtinysp$\nobreakdash-algebra
as $\alpha+u$ (or $\beta+v$, or $\lambda+g$, \dots) without any comment;
in such a~case we will silently assume
that this is the unique $\text{\it scalar}+\text{\it pure\/}$ representation of the element.

\txtskip

Lemma~\ref{lem:scalar-pure-split} says
that the pure part $\pure{A}$ of a~quadratic extension $K\negdtinysp$-algebra $A$
is a~set of representatives of the cosets of the subgroup $K^+$ of the additive group $A^+$.
But more is true\,---\,$\pure{A}$ is a~very special set of representatives.

\thmskip

\begin{lemma}\label{lem:pure-is-subspace}
The pure part of a~quadratic extension\/ $K\negdtinysp$-algebra\/ $A$
is a\/~$K\negdtinysp$-subspace~of\/~$A$.
\end{lemma}

\interskip

\begin{proof}
If $\lambda\in K\negtinysp$ and $u\in\pure{A}$, then $\lambda u\in\pure{A}$.
This is clear if $\lambda=0$ or $u=0\tinysp$;
if $\lambda\neq 0$ and $u\neq 0$,
then $u\notin K\negtinysp$, hence also $\lambda u\notin K\negtinysp$,
and $(\lambda u)^2=\lambda^2u^2\in K\negtinysp$.

Given any $u,\tinysp v\in\pure{A}$, we shall show that $u\narrt+v\in\pure{A}$.

If $u$ and $v$ are linearly dependent,
then $u+v=\lambda u$ for some $\lambda\in K\negtinysp$ or $u+v=\mu v$ for some $\mu\in K\negtinysp$
	(or both),
and in either case $u+v\in\pure{A}$ by what we have proved above.

Now let $u$ and $v$ be linearly independent.

We claim that $1$, $u$, and $v$ are linearly independent;
in particular, $u\narrt+v\notin K\negtinysp$.
Suppose, to the contrary,
that $\alpha+\beta u+\gamma v=0$ for some $\alpha,\tinysp\beta,\tinysp\gamma\in K\negtinysp$,
not all zero.
Then $\alpha\neq 0$ because $u$ and $v$ are linearly independent,
hence at least one of $\beta$, $\gamma$ is non-zero, say $\gamma\neq 0$,
and we have $v=-\alpha/\gamma-(\beta/\gamma)\tinysp u$,
which is a~$\text{\it scalar}+\text{\it pure\/}$ representation of $v\tinysp$;
since the unique such representation of $v$ is $0+v$, we must have $\alpha=0$, a~contradiction.

There exist $\alpha_1,\tinysp\alpha_2\in K\negtinysp$
such that $\beta_1\defeq(u\narrt+v\narrt-\alpha_1)^2\narrt\in K\negtinysp$
and $\beta_{\tinysp2}\defeq(u\narrt-v\narrt-\alpha_2)^2\narrt\in K\negtinysp$,
whence we can express $uv+vu$ in two ways as a~linear combination of $1$, $u$, and $v\tinysp$:
\begin{align*}
uv + vu &\Eq \bigl(\beta_1\narrt-u^2\narrt-v^2\narrt-\alpha_1^2\bigr)
			\widedt+ 2\alpha_1u \widedt+ 2\alpha_1v \\[.5ex]
	&\Eq \bigl(-\beta_{\tinysp2}\narrt+u^2\narrt+v^2\narrt+\alpha_2^2\bigr)
			\widedt- 2\alpha_2\tinysp u \widedt+ 2\alpha_2\tinysp v~.
\end{align*}
Comparing the coefficients at $u$ and $v$
we find that $\alpha_1=-\alpha_2$ and $\alpha_1=\alpha_2=-\alpha_1$,
thus $\alpha_1=0$ and hence $(u\narrt+v)^2=\beta_1\in K\negtinysp$;
since $u\narrt+v\notin K\negtinysp$, we have $u+v\in\pure{A}$.
\end{proof}

\pagebreak[3]\thmskip

We see that every quadratic extension $K\negdtinysp$-algebra $A$ splits as $A=K\negtinysp\oplus\pure{A}$,
where $\pure{A}$ is a~subspace (actually a~hyperplane) of the $K\negdtinysp$-space $A$,
and $\sqr{\pure{A}}\subseteq K\negtinysp$.
This suggests the following definition.

\txtskip

A~\notion{ground Clifford\/ $K\negdtinysp$-algebra}
is an~associative $K\negdtinysp$-algebra of the form $A=K\negtinysp\oplus V\negdtinysp$,
where $V\negdtinysp$~is a~subspace of the $K\negdtinysp$-space~$A$ and $\sqr{V}\subseteq K\negtinysp$.

\txtskip

Recall that a~\notion{Clifford\/ $K\negdtinysp$-algebra} is a~pair $\pair{A,V}$,
where $A$ is a~non-zero unital associative $K\negdtinysp$-algebra,
$V\negdtinysp$~is a~$K\negdtinysp$-subspace of~$A$,
$\sqr{V}\subseteq K1_A=K\negtinysp$, and $V\negdtinysp$~generates~$A\tinysp$
(\ie, $V\negdtinysp$ is a~set of generators of the $K\negdtinysp$-algebra~$A\tinysp$).
More often than not we refer to $\pair{A,V}$ simply as a~Clifford $K\negdtinysp$-algebra~$A$,
knowing that the subspace~$V\negdtinysp$ is somewhere at hand.

If $\pair{A,V}$ is a~Clifford $K\negdtinysp$-algebra,
and $A$ properly contains~$K\negtinysp$ (hence $V\nsubseteq K$),
then $K\negtinysp\inters V=0$.
To prove this, consider any $\alpha\in K\negtinysp\inters V\negdtinysp$,
and pick $v\in V\narrdt\setdiff K\negtinysp$;
since the sum $\alpha+v$ belongs to $V\negdtinysp$,
its square $(\alpha\narrt+v)^2=(\alpha^2\narrt+v^2)+2\tinysp \alpha\tinysp v$ is a~scalar,
thus $\alpha\tinysp v$ is a~scalar, which is possible only if $\alpha=0$.
When $A=K\negtinysp$, the subspace $V\negdtinysp$ is either $0$ or $K\negtinysp$;
in~this special case we always choose $V=0$ instead of $V=K\negtinysp$.

A~Clifford $K\negdtinysp$-algebra $A$
thus always contains the $K\negdtinysp$-subspace $K\negtinysp\oplus V\negdtinysp$,
which is a~foundation on which $A$ is built
as the set of all finite sums of scalars and of elements of~$V\negdtinysp$
and of products of two, three, four, \ldots\ elements of~$V\negtinysp$;
when $V\negdtinysp$ has a~finite dimension~$n$, we can stop with products of~$n$~elements of~$V\negdtinysp$.
A~ground Clifford algebra, as defined above,
is~a~Clifford algebra that consists of nothing else but its foundation;
when we are building such an algebra from its foundation, we never get off the ground.

Let $A=K\negtinysp\oplus V\negdtinysp$ be a~ground Clifford $K\negdtinysp$-algebra.
We denote by $\sqr_{V\negdtinysp}$\label{page:quad-functional-sqrA}
the quadratic functional $V\negdtinysp\to K\negtinysp : v\mapsto v^2$.
Having $\sqr_{V\negdtinysp}$,
we erect the free Clifford $K\negdtinysp$-algebra $\Cliff(\sqr_{V\negdtinysp})$,
which contains $K\negtinysp\oplus V\negdtinysp$ as a~$K\negdtinysp$-subspace.
Elements of $V\negdtinysp$
have the same squares in the algebra $\Cliff(\sqr_{V\negdtinysp})$ as in the algebra~$A$.
However, when $\dim V > 1$,
the algebra $\Cliff(\sqr_{V\negdtinysp})$ properly contains the algebra~$A$,
and the multiplication of elements of $V\negdtinysp$ in $\Cliff(\sqr_{V\negdtinysp})$
differs from their multiplication in~$A$.
There exists a~unique $K\negdtinysp$-algebra homomorphism $h\colon\Cliff(\sqr_{V\negdtinysp})\to A$
which fixes all elements of~$A\tinysp$;
the homomorphism~$h$ collapses the free Clifford algebra $\Cliff(\sqr_{V\negdtinysp})$
onto the ground Clifford algebra~$A$.

\thmskip

\begin{lemma}\label{lem:quad-ext-alg=ground-cliff-alg}
A~quadratic extension\/ $K\negdtinysp$-algebra\/~$A$
is a~ground Clifford\/ $K\negdtinysp$-algebra\/ $K\negtinysp\oplus\pure{A}$.
A~ground Clifford algebra\/ $A=K\negtinysp\oplus V\negdtinysp$ is a~quadratic extension algebra with\/ $\pure{A}=V\negdtinysp$.
\end{lemma}

\interskip

\begin{proof}
The first statement is just a~restatement of Lemma~\ref{lem:pure-is-subspace}.

For the second statement, suppose $A=K\negtinysp\oplus V\negdtinysp$ is a~ground Clifford algebra.
Let $x=\alpha+v$, where $\alpha\in K\negtinysp$ and $v\in V\negdtinysp$, be an arbitrary element of~$A$.
Clearly $x$ satisfies the~quadratic equation $x^2 - 2\tinysp\alpha\tinysp x +(\alpha^2\narrt-v^2)=0$.
Since $x^2=(\alpha^2\narrt+v^2)+2\tinysp\alpha\tinysp v$ is a~scalar
\iff\ $\alpha\tinysp v=0$,
that is,
\iff\ $\alpha=0$ or $v=0$,
it follows that $\pure{A}=V\negdtinysp$.
\end{proof}

\thmskip

Therefore, if~$A$ is a~ground Clifford algebra, it is a~quadratic extension algebra,
the split $A=K\negtinysp\oplus V\negdtinysp$ such that $\sqr V\subseteq K$ is unique,
the subspace $V=\pure{A}$ is determined~by~$A$,
so~we~can write $\sqr_{V\negdtinysp}$ as $\sqr_{\negdtinysp A}$ without ambiguity.
Just keep in mind that the quadratic functional $\sqr_{\negdtinysp A}$
is defined only on the pure part of~$A$, not on the whole~$A$.

\txtskip

How about subalgebras of quadratic extensions algebras?
The following is an immediate consequence of the definitions
of a~quadratic extension algebra and of its~pure~part.

\thmskip

\begin{lemma}\label{lem:subalgs-f-quad-ext-algs}
Let\/ $A$ be a~quadratic extension\/ $K\negdtinysp$-algebra.
Every\/ $K\negdtinysp$-subalgebra\/ $B$ of\/ $A$ is a~quadratic extension\/ $K\negdtinysp$-algebra,
where\/ $\pure{B}=B\inters\pure{A}$.
\end{lemma}

\pagebreak[3]\bigskip

\section{The conjugation and the (reduced) norm}
\label{sec:conjug-and-norm}

\medskip

Let $A$ be a~quadratic extension $K\negdtinysp$-algebra.

\txtskip

Every $x\in A$ has a unique representation $x=\alpha+u$
with $\alpha\in K\negtinysp$ and $u\in\pure{A}\tinysp$;
we shall call $\scl_{\!A}\tinysp x=\scl x\defeq \alpha$
the \notion{scalar part} of~$x$
and $\vct_{\!A}\tinysp x=\vct x\defeq u$
the \notion{pure part} of~$x$.
The mappings $\scl,\tinysp\vct\colon A\mapsto A$ are complementary $K\negdtinysp$-linear projectors
onto the subspaces $K\negtinysp$ resp.~$\pure{A}$ of the $K\negdtinysp$-space~$A\tinysp$:
$\scl A=K\negtinysp$, $\vct A=\pure{A}$,
$\scl^2=\scl$, $\vct^2=\vct$, $\scl\vct=\vct\scl=0$, and $\scl+\vct=\id_A$.

\txtskip

For every $\alpha+u\narrt\in A$ we set $(\alpha\narrt+u)^*\defeq \alpha-u$,
and call the maping $A\narrt\to A \narrt: x\narrt\mapsto x^*$ the~\notion{conjugation} on~$A$.
The conjugation is the endomorphism $\scl-\vct$ of the $K\negdtinysp$-space~$A\tinysp$;
it~is~evidently an involution, that is, $x^{{*}{*}}=x$ for every $x\in A$,
so it is in fact an~auto\-morphism~of the $K\negdtinysp$-space~$A$.
For every $x\in A$ we have
$2\tinysp\scl x=x+x^*$ and $2\tinysp\vct x=x-x^*$,
therefore $x\in K\negtinysp$ iff $\vct x=0$ iff $x^*=x$,
and $x\in\pure{A}$ iff $\scl x=0$ iff $x^*=-x$.

\thmskip

\begin{lemma}\label{lem:pure-u-v-conjug(uv)=vu}
If\/ $u,\tinysp v\in\pure{A}$, then\/ $(uv)^*=vu$.
\end{lemma}

\interskip

\begin{proof}
Let $uv=\alpha+w$, where $\alpha\in K\negtinysp$ and $w\in\pure{A}$.
Since $uv+vu=(u\narrt+v)^2-u^2-v^2$ is a~scalar,
we have $vu=\beta-w$, where $\beta=uv+vu-\alpha\in K\negtinysp$.
We will show that $\beta=\alpha$.

Suppose that $w\neq 0$.
The product $uv\narrdt\cdot vu = u\narrdt\cdot v^2\narrdt\cdot u=u^2\tinysp v^2$ is a~scalar;
because $uv\narrdt\cdot vu = (\alpha\narrt+w)(\beta\narrt-w)
	= (\alpha\beta\narrt-w^2)+(\beta\narrt-\alpha)\tinysp w$,
we must have $(\beta\narrt-\alpha)w=0$, which implies~$\beta-\alpha=0$. 

Now suppose that $w=0$, thus $uv=\alpha$ and $vu=\beta$.
If $u=0$, then $uv=vu=0$, that is, $\alpha=\beta=0$, so suppose that $u\neq 0$.
We have $\alpha u=u\alpha=u\narrdt\cdot uv=u^2 v$ and $\beta u=vu\narrdt\cdot u=vu^2=u^2v$,
thus~$\alpha u=\beta u$, whence $\alpha=\beta$. 
\end{proof}

\thmskip

\begin{lemma}\label{lem:conjug(xy)=conjug(y)conjug(x)}
For all\/ $x,\tinysp y\in A$ we have\/ $(x\tinysp y)^*=y^*x^*$.
\end{lemma}

\interskip

\begin{proof}
Splitting $x=\alpha+u$ and $y=\beta+v$,
we have $(x\tinysp y)^*=(\alpha\beta\narrt+\alpha\tinysp v\narrt+\beta u\narrt+uv)^*
	= (\alpha\beta)^*+(\alpha\tinysp v)^*+(\beta u)^*+(uv)^*
	= \alpha\beta-\alpha\tinysp v-\beta u+vu
	= y^*\negtinysp x^*$.
\end{proof}

\pagebreak[3]
\thmskip

If $x,\tinysp y\in A$ and $\alpha\defeq x\narrt+y\in K\negtinysp$, then $x$ and $y=\alpha-x$ commute;
in particular, $x$~and~$x^*$ always commute.
For every $x\in A$ we call $\Qnorm_{\!A}(x)=\Qnorm(x)\defeq x\tinysp x^*=x^*\negtinysp x\in K\negtinysp$
the~\notion{(reduced) norm of\/~$x$}.
If~$x=\alpha+u$, then $\Qnorm(x)=\alpha^2-u^2$.
For a~pure $u$ we have $\Qnorm(u)=-u^2$,
that is,
denoting the restriction of $\Qnorm_{\!A}$ to $\pure{A}$ by $\qnorm_{\!A}$,
we have $\qnorm_{\!A} = -\sqr_{\negdtinysp A}$.

\thmskip

\begin{lemma}\label{lem:norm-is-multiplicative}
The norm\/ $\Qnorm_{\!A}\colon A^{\tinysp\cdot}\to K^{\tinysp\cdot}\negtinysp$
is a~homomorphism of multiplicative monoids.
\end{lemma}

\interskip

\begin{proof}
Clearly $\Qnorm(1)=1$.
If $x,\tinysp y\narrt\in A$,
then $\Qnorm(x\tinysp y) = (x\tinysp y)(x\tinysp y)^*
	= x\tinysp y\tinysp y^*\negtinysp x^*
	= x\narrdt\cdot\tinysp y\tinysp y^*\negdtinysp\narrdt\cdot x^*
	= x\tinysp x^*\negdtinysp\narrdt\cdot y\tinysp y^*
	=\Qnorm(x)\tinysp\Qnorm(y)$.
\end{proof}

\pagebreak[3]
\thmskip

\begin{lemma}\label{lem:invertibility-in-quatalg}
An element\/ $x$ of\/ $A$ is invertible \iff\/ $\Qnorm(x)\neq 0$.
\end{lemma}

\interskip
\begin{proof}
If $x$ is invertible,
then $\Qnorm(x)\tinysp\Qnorm(x^{-1})=\Qnorm(x\tinysp x^{-1})=\Qnorm(1)=1$,
hence $N(x)\neq 0$.
Conversely, if~$\Qnorm(x)\neq\nolinebreak 0$,
then $x$ is invertible with the inverse $x^*\negdtinysp/\Qnorm(x)$.
\end{proof}

\pagebreak[3]\txtskip

The following proposition summarizes some properties of conjugation.

\thmskip

\begin{proposition}\label{prop:conjug-is-antiautomorph}
The conjugation\/ $x\mapsto x^*$ of a~quadratic extension algebra\/ $A$
is an involutive anti-automorphism of\/~$A$ which fixes every scalar.
Moreover, $x+x^*$ and\/ $x\tinysp x^*$ are scalars for every\/ $x\in A$.
\end{proposition}

\thmskip

We can characterize quadratic extension algebras by existence of a~`conjugation'.
First, a preparatory lemma.

\thmskip

\begin{lemma}\label{lem:(x+xpr)-and-(x*xpr)-scalars}
Let\/ $A$ be a~quadratic extension algebra.
If\/ $x,\tinysp\xpr\in A$
and\/ $x$ is not a~scalar, while\/ $x+\xpr$ and\/ $x\tinysp\xpr$ are scalars,
then\/ $\xpr=x^*\negtinysp$.
\end{lemma}

\interskip

\begin{proof}
Splitting $x=\alpha+u$, $\xpr=\beta-u$,
we have a~scalar $x\tinysp\xpr=(\alpha\beta\narrt-u^2)+(\beta\narrt-\alpha)\tinysp u$,
thus $\beta-\alpha=0$ because $u\neq 0$.
\end{proof}

\thmskip

Note that if $x,\tinysp\xpr\in A$ are such that $x+\xpr$ is a~scalar,
and $x$ is a~scalar, then $\xpr$ is a~scalar.
Also, given a~scalar $x$, both $x+\xpr$ and $x\tinysp\xpr$ are scalars for \emph{any} scalar~$\xpr$.

And now, here is the promised characterization.

\thmskip

\begin{proposition}\label{prop:quadext-algs-by-conjug}
Let\/ $A$ be a~non-zero unital associative\/ $K\negdtinysp$-algebra,
and let\/ $A\to A : x\mapsto\overbar{x}$ be a~function
such that\/ $x+\overbar{x}$ and\/ $x\tinysp\overbar{x}$ are scalars for every\/ $x\in A$.
Then\/~$A$ is~a~quadratic extension\/ $K\negdtinysp$-algebra,
and\/ $\overbar{x}=x^*$ for every\/ $x\narrt\in A\narrdt\setdiff K\negtinysp$.
Consequently, if\/~$\overbar{\alpha}=\alpha$ for every\/ $\alpha\in K$,
then\/ $\overbar{x}=x^*$ for every\/ $x\in A$.
\end{proposition}

\interskip

\begin{proof}
Every $x\in A$ satisfies the equation $x^2-x\tinysp(x\narrt+\overbar{x})+x\tinysp\overbar{x}=0$,
which proves that $A$~is~a~quadratic extension $K\negdtinysp$-algebra.
If~$x\in A\narrdt\setdiff K\negtinysp$, then $\overbar{x}=x^*$
by~Lemma~\ref{lem:(x+xpr)-and-(x*xpr)-scalars}.
\end{proof}

\pagebreak

\section{Examples of quadratic extension algebras}
\label{sec:examps-quadext-algs}

\medskip

Every quaternion algebra $\tQalgLam{a}{b}{K}$, where $a,\tinysp b\in\punct{K\tinysp}$,
is a~quadratic extension $K\negdtinysp$-algebra.

\txtskip

The field $K\negtinysp$ itself is a~quadratic extension $K\negdtinysp$-algebra
whose pure part is~$0$.

\txtskip

The quadratic extension $K\negdtinysp$-algebras $A$
with a~one$\tinysp$-dimensional $\pure{A}=K\negtinysp u$
are precisely the free Clifford algebras generated by one$\tinysp$-dimensional quadratic spaces.
The kind of algebra we get depends on the nature of the scalar~$u^2$.
If $u^2=0$, then $A$ is 
the algebra of `dual scalars' over the field~$K\negtinysp$.
If $u^2$ is a~square of a~non-zero scalar,
then $A$ is isomorphic to the $K\negdtinysp$-algebra $K\negtinysp\narrdt\times K\negtinysp$.
Finally, if $u^2=\alpha$ is not a~square of a~scalar,
then $A$ is (isomorphic to)
the quadratic extension field $K\negtinysp\bigl(\sqrt{\alpha\tinysp}\dtinysp\bigr)$ of~$K\negtinysp$.

\txtskip

Next on our agenda
are the quadratic extension $K\negdtinysp$-algebras $A$ with a~two$\tinysp$-dimensional pure part;
we shall determine all of them 
and classify them up to an isomorphism.

Let $\pair{e,f}$ be an orthogonal basis of the quadratic space $\bigpair{\pure{A},\sqr_{\negdtinysp A}}$.
The~symmetric bilinear functional
associated with the quadratic functional~$\sqr_{\negdtinysp A}\colon\pure{A}\to K : u\mapsto u^2$
is
\begin{equation*}
B_{\negdtinysp A}(u,v)
	\Defeq \thalf\bigl((u\narrt+v)^2-u^2-v^2\bigr)
	\Eq \thalf\bigl(uv+vu\bigr)
	\Eq \thalf\bigl(uv+(uv)^*\bigr)
	\Eq \scl(uv)~.
\end{equation*}
Since $e\negtinysp f+f\negtinysp e = 2\tinysp B_{\negdtinysp A}(e,f) = 0$,
we see that~$\scl(e\negtinysp f)=0$, which means that $e\negtinysp f$~is~pure.
We have $e^2=\alpha$, $f^2=\beta$, $e\negtinysp f=-\negtinysp f\negtinysp e=\gamma\tinysp e+\delta f$
for some scalars $\alpha$, $\beta$, $\gamma$, $\delta$.
The necessary and sufficient conditions
$(e\tinysp e)f=e(e\negtinysp f)$, $(e\negtinysp f)e=e(f\negtinysp e)$,~\ldots\
for associativity of~$A$ give us the equations
$\alpha=\delta^2$, $\beta=\gamma^2$, and $\alpha\tinysp\gamma=\beta\tinysp\delta=\gamma\tinysp\delta=0$.

When $\alpha\neq 0$, we have $\delta\neq 0$ and $\beta=\gamma=0$
(the case $\beta\neq0 \Implies \gamma\neq 0 \Implies \alpha=\delta=0$ is analogous).
Renaming $e/\delta$ as $e$, we obtain a~basis $\pair{e,f}$ of~$\pure{A}$
such that $e^2=1$, $f^2=0$, and $e\negtinysp f=-\negtinysp f\negtinysp e=f$.
Let the matrices $I=\bigtxtmtx{\begin{smallmatrix}1 & 0\\[.2ex]0 & 1\end{smallmatrix}}$,%
~$E=\bigtxtmtx{\begin{smallmatrix}1 & 0\\[.2ex]0 & -1\end{smallmatrix}}$,%
~and~$F=\nolinebreak\bigtxtmtx{\begin{smallmatrix}0 & 1\\[.2ex]0 & 0\end{smallmatrix}}$
have entries in~$K\negtinysp$;
they form a~basis of the underlying $K\negdtinysp$-space of the $K\negtinysp$\nobreakdash-algebra
$R=\bigtxtmtx{\begin{smallmatrix}\!K & \!K\!\\[.2ex]\!0 & \!K\!\end{smallmatrix}}$
of~the~upper triangular $2\narrdt\times 2$ matrices with entries in~$K\negtinysp$.
The~matrices $E$~and~$F$ satisfy the equations $E^2=I$, $F^2=0$, and $EF=-FE=F$,
thus the $K\negtinysp$\nobreakdash-linear mapping $A\to R$
which sends $1$,~$e$,~$f$ respectively to $I$, $E$, $F$
is an isomorphism of $K\negtinysp$\nobreakdash-algebras.
Every~matrix~$X\negtinysp=\nolinebreak\negtinysp
	\bigtxtmtx{\begin{smallmatrix}\alpha & \beta\\[.2ex]0 & \gamma\end{smallmatrix}}\narrt\in R$
has the scalar part $\scl X \narrt= \thalf(\alpha\narrt+\gamma) \narrt= \thalf\tr X$
{\large(}actually $\scl X \negtinysp =
	\nolinebreak\negtinysp \bigl(\thalf\tr X\bigr)\narrdt\cdot I\dtinysp${\large)},
the norm $\Qnorm(X) \narrt= \alpha\tinysp\gamma \narrt= \det X$,
and the conjugate
$X^* \narrt= \bigtxtmtx{\begin{smallmatrix}\gamma & -\beta\\[.2ex]0 & \alpha\end{smallmatrix}}
	\narrt= \adjug{X}$
{\large(}where for any positive integer $n$ and any matrix $M\in\Malg_n(K)$,
$\adjugM$ denotes the adjugate of the matrix~$M$,
that is, the transpose of the matrix of cofactors of~$M\tinysp${\large)}.
The pure part of $R$ is $\pure{R}=\set{X\narrt\in R\suchthat\tr X=0}$,
that is, $\pure{R}$ consists of all matrices of the form
$\bigtxtmtx{\begin{smallmatrix}\alpha & \beta\\[.2ex] 0 & -\alpha\end{smallmatrix}}
	= \alpha\tinysp E+\beta F$
with $\alpha,\tinysp\beta\in K$.
The quadratic functional $\sqr_{\negdtinysp A}$
is isomorphic to the quadratic functional $p_2\colon K^2\to K:\tuple{x,y}\mapsto x^2$,
and the free Clifford algebra $\Cliff(p_2)$%
\,---\,isomorphic to the $K\negdtinysp$-algebra $\tQalgLam{1}{0}{K}$\,---\,%
has a basis $\tuple{1,i,j,k}$ with the multiplication table
$i^2=1$, $ij=-ji=k$, $ik=-ki=j$, $j^2=k^2=jk=kj=0$.
The `collapsing'~homo\-morphism of $K\negdtinysp$-algebras $h\colon\Cliff(p_2)\to A$,
which sends $i$ to $e$ and $j$ to~$f$,
sends $k=ij$ to~$e\negtinysp f=f$,
thus the kernel of $h$ is the one-dimensional subspace $K\negdtinysp\narrt\cdot(k\narrt-j)$
of $\Cliff(p_2)$.

\pagebreak[3]

The only other case is $\alpha=\beta=0$, when also $\delta=\gamma=0\tinysp$;
now $\pure{A}$ is a~zero$\tinysp$-product algebra,
and adjoining $1$ to it we obtain the algebra~$A$.


We have found that there are only two isomorphism classes
of three$\tinysp$-dimensional quadratic extension $K\negdtinysp$-algebras:
one class is of the $K\negdtinysp$-algebra
obtained by adjoining the unity to the zero$\tinysp$-product algebra on $K^2\negtinysp$,
and the other class is of the $K\negdtinysp$-algebra
of~the~upper triangular $2\narrdt\times 2$ matrices with entries in~$K\negtinysp$.

\txtskip

Digressing for a~couple of moments,
let us consider an arbitrary three-dimensional $K\negtinysp$\nobreakdash-algebra~$A$.
We~have just discussed the case where every element of $A$ satisfies a~quadratic equation over~$K$.
Otherwise, there is an element $a$ of $A$ that does not satisfy any quadratic equation.
Since $A$ is three-dimensional, $a$~does satisfy some cubic equation,
thus $f(a)=0$ for some monic polynomial $f(X)\in K[X]$ of degree~$3$.
Moreover, every polynomial $g(X)\in K[X]$ such that $g(a)=0$ is divisible by~$f(X)$.
The elements $1$, $a$, $a^2$ of $A$ are linearly independent,
and $A = K\narrt\oplus Ka\narrt\oplus Ka^2 \isomorph K[X]/f(X)\tinysp K[X]$.
We know how the structure of $A$ can be worked out in detail,
depending on factorization of $f(X)$ into irreducible factors;
we won't go into this, we just point out that $A$ is commutative.
So~we have this nice little result:
a~three-dimensional $K\negdtinysp$-algebra that is not commutative
is isomorphic to the $K\negdtinysp$-algebra
of upper-triangular $2\narrdt\times2$ $K\negdtinysp$-matrices.
Since the latter $K\negdtinysp$-algebra possesses divisors of zero, we have also the following corollary:
if~a~three-dimensional $K\negdtinysp$-algebra is a~domain, then it is a~field.%
\footnote{Every finite-dimensional algebra is algebraic,
and every algebraic algebra that is a~domain (has no divisors of zero) is a~division ring.}

\txtskip

As the~last example we provide a~generous supply of rather degenerate quadratic extension algebras:
if we adjoin $1$ to any zero$\tinysp$-square $K\negdtinysp$-algebra~$B$,
we obtain a~quadratic extension $K\negdtinysp$-algebra $A=K\negtinysp\oplus B$
whose pure part is~$B$;
in particular, we may choose $B$ to be any zero$\tinysp$-product $K\negdtinysp$-algebra.

There do exist zero$\tinysp$-square algebras that are not zero$\tinysp$-product algebras.
For example, the $K\negdtinysp$-algebra $B\defeq Ke\oplus Kf\oplus Kg$ with the multiplication table
$e^2=f^2=g^2=0$ and
$e\negtinysp f=-\negtinysp f\negtinysp e=f\negtinysp g=-g\negtinysp f=g\tinysp e=-eg=e+f+g$
is a~zero$\tinysp$-square algebra.
Choosing the basis $\tuple{e,f\negtinysp,e\negtinysp f}$ of~$B$
instead of the basis $\tuple{e,f\negtinysp,g}$,
we discover that the quadratic extension $K\negtinysp$\nobreakdash-algebra $A\defeq K\negtinysp\oplus B$
is isomorphic to the exterior algebra $\ExtAlg\bigl(K^{\tinysp2}\bigr)$;
put slightly differently,
$A$ is isomorphic to the free Clifford algebra $\Cliff(o_2)$
of the zero quadratic functional $o_2$ on~$K^{\tinysp2}$,
where $\Cliff(o_2)$ can be also presented as~$\tQalgLam{0}{0}{K}$.

The variety of zero$\tinysp$-square $K\negdtinysp$-algebras
is closed under taking direct products,~subalge\-bras,~and homomorphic images.
The not\nobreakdash-completely-trivial zero\nobreakdash-square
	$K\negdtinysp$-algebra $B$
generates a~subvariety of this variety
which contains a~great many not\nobreakdash-com\-pletely-trivial
	zero$\tinysp$-square $K\negdtinysp$-algebras.%
\footnote{Question: does $B$ generate the entire variety of zero-square $K\negdtinysp$-algebras?}
%

\pagebreak

\section{Pure calculus on a~quadratic extension algebra}
\label{sec:pure-calculus}

\medskip

When I first came acrosss Hamilton's quaternions,
they appeared to me as a~not entirely natural marriage
of the real line and the oriented three-dimensional Euclidean space.
Every quaternion was a~sum of a~scalar and a~vector,
and the product of two such sums $s+x$ and $t+y$ was given by the formula
\begin{equation}\label{eq:hamilton-quaterns-product}
(s\narrt+x)(t\narrt+y)
	\Eq \bigl(s\tinysp t - \bilin{x}{y}\bigr)
		\wide+ \bigl(tx + sy + x\narrdt\cross y\bigr)
\end{equation}
which was a concoction of several ingredients of the algebra of three-dimensional vectors:
a product of two scalars, two products of a~vector by a~scalar,
a scalar/dot product of vectors, a~vector/cross products of vectors,
a~sum of two scalars, and a~sum of three vectors.
For quite a~while I~thought that this was how the quaternions had actually been designed.
I~was mistaken, of course.
In~the beginning there were quaternions, quadruples of real numbers,
created  by Hamilton on October~16 of the year 1843,
as~a~generalization of complex numbers
(after almost 13~years\,---\,some say that there were `only' 8~years\,---\,%
of~futile attempts to set up such a~generalization with triples of real numbers).
When~the quaternions were $40$ years old, give or take a~year or two,
they were mercilessly hacked into their scalar and vector parts,
and then these chunks (still dripping blood)
were used to build
the~vector calculus of the oriented three-dimensional Euclidean space:
\begin{quote}
It is probably true that Hamilton spent too much time on quaternions.
He~did~little else until his death in 1865, and few mathematicians shared his enthusiasm.
Nevertheless quaternions changed the course of mathematics, though not in the way Hamilton intended.
In the 1880s Josiah Willard Gibbs and Oliver Heaviside created what we now know as vector analysis,
essentially by separating the real (``scalar'') part of a~quaternion
from its imaginary (``vector'') part.
Hamilton's followers were outraged to see the simple and elegant quaternions torn limb from limb,
but the idea caught on with physicists and engineers, and it still holds~sway~today.

[John Stillwell \cite{StillwellMaIH}, pp.~402--403.]
\end{quote}
We are now going to follow in steps of Gibbs and Heaviside,
constructing an analogue of (the~algebraic part of) the~``vector calculus''
from the scalar and pure parts of elements of an arbitrary quadratic extension algebra;
this will be the ``pure calculus'' of this section's title.
To~get the right idea where and how to start,
we take a look at the product of pure Hamilton's quaternions (that is, vectors) $x$ and $y$,
which is
\begin{equation*}
x\tinysp y \Eq -\negdtinysp\bilin{x}{y}\widedt+x\narrdt\cross y~,
\end{equation*}
where the first summand is the scalar part, and the second summand is the vector part,
of the product.

\pagebreak[3]\txtskip

Let $A$ be a~quadratic extension $K\negdtinysp$-algebra.

\txtskip

For any two pure elements $u$ and $v$ of~$A$ we define
%
\begin{align*}
\bilin{u}{v}
	&\Defeq -\scl(uv)
		\Eq -\thalf(uv\narrt+vu)
		\Eq -B_{\negdtinysp A}(u,v) \wide\in K~, \\[.5ex]
u\narrdt\cross v
	&\Defeq \vct(uv)
		\Eq \thalf(uv\narrt-vu) \wide\in \pure{A}~.
\end{align*}
We should perhaps write the two operations as
$\bilin{u}{v}_{\!A}$ and $u\narrt{\cross_{\!A}}v$
to avoid confusion\,---\,but we will always know what we are doing, won't we?
Using the \notion{scalar product} $\bilin{\anon}{\anon}$
and the \notion{pure product} $\anon\narrdt\cross\anon$ of pure elements,
the product of any two elements $\alpha+u$ and $\beta+v$ of~$A$
can be written in the $\text{\it scalar}+\text{\it pure\/}$ form as
\begin{equation}\label{eq:general-quadextalg-product}
(\alpha\narrt+u)(\beta\narrt+v)
	\Eq \bigl(\alpha\beta - \bilin{u}{v}\bigr)
		\wide+ \bigl(\beta u + \alpha\tinysp v + u\narrdt\cross v\bigr)~,
\end{equation}
which is (of course)
precisely the no-longer-so-strange-looking formula~\eqref{eq:hamilton-quaterns-product}.
Note that the scalar product is a~bilinear functional on~$\pure{A}$,
and that the pure multiplication is a~bilinear operation on~$\pure{A}\tinysp$;
moreover, $\bilin{\anon}{\anon}$ is symmetric, $\bilin{v}{u}=\bilin{u}{v}$,
while $\anon\narrdt\cross\negtinysp\anon$ is anticommutative,
$v\narrdt\cross u=-\tinysp u\narrdt\cross v$.
Note that
\begin{equation*}
\bilin{u}{u} \Eq -u^2 \Eq \Qnorm(u)
\end{equation*}
for every pure~$u$,
thus $\bilin{\anon}{\anon}$ is the symmetric bilinear functional
associated with the quadratic functional $\qnorm_{\!A}=-\sqr_{\negdtinysp A}$ on the subspace~$\pure{A}$
(which is the restriction of the quadratic functional $\Qnorm_{\!A}$ defined on the space~$A$).

\txtskip

Let $u,\tinysp v,\tinysp w\in\pure{A}$.

\txtskip

We say that $u$ and $v$ are \notion{orthogonal}, and write $u\perp v$, if $\bilin{u}{v}=0$.
Since $\bilin{u}{v}=-\scl(uv)=-\scl(vu)$, we have the following equivalences:
\begin{equation*}
u\perp v 
	\wide\Isequiv \text{$uv$ is pure}
	\wide\Isequiv \text{$vu$ is pure}
	\wide\Isequiv uv = u\narrdt\cross v~.
\end{equation*}
Since $\bilin{u}{v}=-\thalf(uv\narrt+vu)$,
$u$ and $v$ are orthogonal \iff\ they anticommute, that is, \iff\ $vu=-uv$.
Similarly, since $u\narrdt\cross v=\thalf(uv\narrt-vu)$,
we have $u\narrdt\cross v=0$ \iff\ $u$ and $v$ commute.

\txtskip

Clearly $u\narrdt\cross u=0$.
Slightly more generally,
if $u$ and $v$ are linearly dependent,
then $u=\lambda v$ for some scalar $\lambda$
or $v=\mu\tinysp u$ for some scalar $\mu$ (or both),
and in either case $u\narrdt\cross v=0$.
Equivalently, if $u\narrdt\cross v\neq 0$, then $u$ and $v$ are linearly independent.

The pure product is half the additive commutator in an~associative algebra,
thus $\bigpair{\pure{A},{\cross}}$ is~a~Lie algebra,
hence $\anon\narrdt\cross\negtinysp\anon$ satisfies the Jacobi's identities:
\begin{align*}
u\narrdt\cross(v\narrdt\cross w)
	+ v\narrdt\cross(w\narrdt\cross u)
	+ w\narrdt\cross(u\narrdt\cross v)
	& \Eq 0~, \\[1ex]
(u\narrdt\cross v)\narrdt\cross w
	+ (v\narrdt\cross w)\narrdt\cross u
	+ (w\narrdt\cross u)\narrdt\cross v
	& \Eq 0~.
\end{align*}
%

\pagebreak[3]\txtskip

Here is another easily obtainable identity:
on the one hand
	$\Qnorm(uv)=\Qnorm(u)\tinysp\Qnorm(v)=(-u^2)(-v^2)=u^2\tinysp v^2$,
and on the other hand
	$\Qnorm(uv)=\bigl(\scl(uv)\bigr)^{\negtinysp2}-\bigl(\vct(uv)\bigr)^{\negtinysp2}\negtinysp$,
therefore
\begin{equation}\label{eq:u^2v^2=<u,v>^2-(uxv)^2}
u^2\tinysp v^2 \Eq \bilin{u}{v}^{\negtinysp2} - (u\narrdt\cross v)^2~.
\end{equation}

\vspace{-1ex}
\txtskip

For a~few moments assume that\/ $V\negdtinysp$ is some vector space over~$K$,
that $\bilin{\anon}{\anon}$ is a~$K\negdtinysp$\nobreakdash-bilinear functional on~$V\negdtinysp$
(not necessarily symmetric),
and that $\anon\narrdt\cross\negtinysp\anon$ is a~$K\negdtinysp$-bilinear operation on $V\negdtinysp$
(not necessarily anticommutative).
Define the bilinear multiplication $\pair{x,y}\mapsto x\tinysp y$ on $K\negtinysp\oplus V\negdtinysp$
by%
\footnote{This is in fact the most general situation.
Meaning what?
Let $\pair{x,y}\mapsto x\tinysp y$ be any bilinear operation on $K\negtinysp\oplus V\negdtinysp$
with the neutral element (multiplicative identity)~$1\oplus 0$.
Denoting by $\kappa\colon K\negtinysp\oplus V\to K$ and $\varphi\colon K\negtinysp\oplus V\to V\negdtinysp$
the projections onto the direct summands,
we set $\bilin{u}{v}\defeq -\kappa(uv)$ and $u\narrdt\cross v\defeq\varphi(uv)$
for any $u,\tinysp v\in V\negdtinysp$;
then~the product $(\alpha\oplus u)(\beta\oplus v)$
is computed according to formula~\eqref{eq:general-bilinear-multiplication}.}
\begin{equation}\label{eq:general-bilinear-multiplication}
(\alpha\oplus u)(\beta\oplus v)
	\Defeq \bigl(\alpha\beta-\bilin{u}{v}\bigr)
			\wide\oplus \bigl(\beta u+\alpha v+u\narrdt\cross v\bigr)~.\end{equation}
Let us compute the associator of\/
$x_i=\alpha_i\oplus u_i \in K\oplus V\negdtinysp$, $\dtinysp i=1,\,2,\,3\tinysp$:
\begin{align*}
(x_1 x_2)\tinysp x_3-x_1(x_2\tinysp x_3)
	&= \bigl(-\negdtinysp \bilin{u_1\narr\cross u_2}{\tinysp u_3}
			+ \bilin{u_1}{\tinysp u_2\narr\cross u_3}\bigr) \\[.5ex]
	&\hphantom{{}={}}
		~~\oplus \bigl((u_1\narr\cross u_2)\narr\cross u_3
				\narrt- \bilin{u_1}{\negtinysp u_2}\negtinysp u_3
				\narrt- u_1\narr\cross(u_2\narr\cross u_3)
				\narrt+ \bilin{u_2}{\negtinysp u_3}\negtinysp u_1\bigr)\,.
\end{align*}
We see that the multiplication defined on $K\oplus V\negdtinysp$ is associative \iff\
\begin{gather*}
\bilin{u_1\narr\cross u_2}{\tinysp u_3}
	\Eq \bilin{u_1}{\tinysp u_2\narr\cross u_3}~, \\[1ex]
(u_1\narr\cross u_2)\narr\cross u_3 - u_1\narr\cross(u_2\narr\cross u_3)
	\Eq \bilin{u_1}{\negtinysp u_2}\negtinysp u_3 - \bilin{u_2}{\negtinysp u_3}\negtinysp u_1
\end{gather*}
for all\/ $u_1,\tinysp u_2,\tinysp u_3\in V\negdtinysp$.
Assuming that $\pair{V,{\cross}}$ is a~Lie algebra,
we can rewrite the left hand side of the second identity above as
\begin{align*}
(u_1\narrdt\cross u_2)\narrdt\cross u_3 - u_1\narrdt\cross(u_2\narrdt\cross u_3)
	&\Eq (u_1\narrdt\cross u_2)\narrdt\cross u_3 + (u_2\narrdt\cross u_3)\narrdt\cross u_1 \\[.5ex]
	&\Eq -\tinysp(u_3\narrdt\cross u_1)\narrdt\cross u_2~.
\end{align*}
In our case the multiplication on $A=K\oplus\pure{A}$ is associative
and $\bigpair{\pure{A},{\cross}}$ is a~Lie algebra,
so the following identities hold for all $u,\tinysp v,\tinysp w\in\pure{A}\tinysp$:
\begin{gather}
\bilin{u\narr\cross v}{\tinysp w} \Eq \bilin{u}{\tinysp v\narr\cross w}~,
		\label{eq:<uxv,w>=<u,vxw>} \\[1ex]
(u\narrdt\cross v)\narrdt\cross w \Eq \bilin{u}{w}v - \bilin{v}{w}u~,
		\label{eq:(uxv)xw=<u,w>v-<v,w>u} \\[1ex]
u\narrdt\cross(v\narrdt\cross w) \Eq \bilin{u}{w}v - \bilin{u}{v}w~;
		\label{eq:ux(vxw)=<u,w>v-<u,v>w}
\end{gather}
identity~\eqref{eq:ux(vxw)=<u,w>v-<u,v>w} follows from identity~\eqref{eq:(uxv)xw=<u,w>v-<v,w>u}
because
$u\narrdt\cross(v\narrdt\cross w)
	= -(-w\narrdt\cross v)\narrdt\cross u
	= (w\narrdt\cross v)\narrdt\cross u$
and because $\bilin{\anon}{\anon}$ is symmetric.
Using the identities~\eqref{eq:<uxv,w>=<u,vxw>} and~\eqref{eq:(uxv)xw=<u,w>v-<v,w>u},
we can derive the identity
\begin{equation}\label{eq:<u1xv1,u2xv2>=...}
\bilin{u_1\narrdt\cross u_2}{\dtinysp v_1\narrdt\cross v_2}
	\Eq	\begin{vmatrix}
		\,\bilin{u_1}{v_1} & \bilin{u_1}{v_2}\, \\[.5ex]
		\,\bilin{u_2}{v_1} & \bilin{u_2}{v_2}\,
		\end{vmatrix}%
		~,
\end{equation}
which holds for all $u_1,\tinysp u_2,\tinysp v_1,\tinysp v_2\in\pure{A}\tinysp$:
\begin{align*}
\bilin{u_1\narrdt\cross u_2}{\dtinysp v_1\narrdt\cross v_2}
	&\Eq \bigbilin{(u_1\narrdt\cross u_2)\narrdt\cross v_1}{\dtinysp v_2} \\[.5ex]
	&\Eq \bigbilin{\bilin{u_1}{v_1}\negtinysp u_2\narrt-\bilin{u_2}{v_1}\negtinysp u_1}%
									{\dtinysp v_2} \\[.5ex]
	&\Eq \bilin{u_1}{v_1}\bilin{u_2}{v_2} \widedt- \bilin{u_2}{v_1}\bilin{u_1}{v_2}~.
\end{align*}
Taking $u_1\defeq v_1\defeq u$ and $u_2\defeq v_2\defeq v$ in the identity~\eqref{eq:<u1xv1,u2xv2>=...},
we obtain the identity~\eqref{eq:u^2v^2=<u,v>^2-(uxv)^2}.

\pagebreak[3]\txtskip

We define the \notion{mixed product},
a trilinear functional $\trilin{\anon}{\anon}{\anon}$ on~$\pure{A}$,
by
\begin{equation*}
\trilin{u}{v}{w}
	\Defeq \bilin{u\narrdt\cross v}{\tinysp w}
		\Eq \bilin{u}{\tinysp v\narrdt\cross w}~,
\end{equation*}
for all $u,\tinysp v,\tinysp w\in\pure{A}$.
Since clearly $\trilin{u}{u}{v}=\trilin{u}{v}{v}=0$ for all $u,\tinysp v\in\pure{A}$,
the mixed product is an \emph{alternating} trilinear functional,
thus also $\trilin{u}{v}{u}=0$ for all $u,\tinysp v\in\pure{A}$,
and
\begin{equation*}
\trilin{u_{\sigma1}}{\tinysp u_{\sigma2}}{\tinysp u_{\sigma3}}
	\Eq \sgn(\sigma)\narrdt\cdot\trilin{u_1}{u_2}{u_3}
\end{equation*}
for all $u_1,\tinysp u_2,\tinysp u_3\in\pure{A}$ and every permutation $\sigma$ of $\set{1,2,3}$.
If $u,\tinysp v,\tinysp w\in\pure{A}$ are linearly dependent, then $\trilin{u}{v}{w}=0\tinysp$;
indeed, if, say, $w=\alpha\halftinysp u+\beta v$ for some $\alpha,\tinysp\beta\in K$,
then
\begin{equation*}
\trilin{u}{v}{w}
	\Eq \trilin{u}{v}{\tinysp\alpha\halftinysp u\narrdt+\beta v}
	\Eq \alpha\narrdt\cdot\trilin{u}{v}{u} + \beta\narrdt\cdot\trilin{u}{v}{v}
	\Eq 0.
\end{equation*}
Equivalently, if $\trilin{u}{v}{w}\neq 0$, then $u$, $v$, $w$ are linearly independent.

For all $u,\tinysp v\in\pure{A}$ we have
$\bilin{u}{u\narrdt\cross v} = \trilin{u}{u}{v} = 0$
and also $\bilin{v}{u\narrdt\cross v} = 0$,
that is,
\begin{equation*}
u\perp u\narrdt\cross v~, \qquad v\perp u\narrdt\cross v~.
\end{equation*}

We conclude the section with the identity
\begin{equation*}
\trilin{u\narrdt\cross v}{\dtinysp v\narrdt\cross w}{\dtinysp w\narrdt\cross u}
	\Eq \trilin{u}{v}{w}^2~,
\end{equation*}
which holds for all pure $u$, $v$, and~$w\tinysp$:
\begin{align*}
\trilin{u\narrdt\cross v}{\dtinysp v\narrdt\cross w}{\dtinysp w\narrdt\cross u}
	&\Eq \bigbilin{(u\narrdt\cross v)\narrdt\cross(v\narrdt\cross w)}{\dtinysp w\narrdt\cross u}
												\\[.5ex]
	&\Eq \bigbilin{\bilin{u}{\tinysp v\narrdt\cross w}\negtinysp v
			- \bilin{v}{\tinysp v\narrdt\cross w}\negtinysp u}{\dtinysp w\narrdt\cross u}
												\\[.5ex]
	&\Eq \bigbilin{\trilin{u}{v}{w}\negtinysp v}{\dtinysp w\narrdt\cross u} \\[.5ex]
	&\Eq \trilin{u}{v}{w}\narrdt\cdot\bilin{v}{\tinysp w\narrdt\cross u} \\[.5ex]
	&\Eq \trilin{u}{v}{w}\narrdt\cdot\trilin{v}{w}{u} \\[.5ex]
	&\Eq \trilin{u}{v}{w}^2~.
\end{align*}

\vspace{-1.5ex}\pagebreak[3]\bigskip

\section{The `grand identity'}
\label{sec:grand-identity}

\medskip

Let $A$ be a~quadratic extension $K\negdtinysp$-algebra.

\txtskip

Our aim in this section is to prove the identity
\begin{equation}\label{eq:grand-identity}
\begin{vmatrix}
\,\bilin{u_1}{v_1} & \bilin{u_1}{v_2} & \bilin{u_1}{v_3}\, \\[.5ex]
\,\bilin{u_2}{v_1} & \bilin{u_2}{v_2} & \bilin{u_2}{v_3}\, \\[.5ex]
\,\bilin{u_3}{v_1} & \bilin{u_3}{v_2} & \bilin{u_3}{v_3}
\end{vmatrix}
	\wide\Eq
\trilin{u_1}{u_2}{u_3}\narrdt\cdot\trilin{v_1}{v_2}{v_3}~,
\end{equation}
which holds for all\/ $u_i,\tinysp v_i\in\pure{A}$, $\dtinysp i=1,\, 2,\, 3$.%
\footnote{Suppose that $\pure{A}=K^3\negtinysp$,
that $\bilin{u}{v}$ is the usual scalar product
and $u\narrdt\cross v$ is the usual vector product of triples $u,\tinysp v\narrt\in K^3$,
so that $\trilin{u}{v}{w}$ is the determinant of the matrix with the rows (or columns) $u$,~$v$,~$w$.
Denote by $U\negtinysp$ the matrix with the rows $u_1$, $u_2$, $u_3$
and by $V\negdtinysp$ the matrix with the columns $v_1$,~$v_2$,~$v_3$;
then the identity~\eqref{eq:grand-identity} says that $\det(U\tinysp V)=\det(U)\det(V)$.}

\txtskip

This time we cannot prove the identity by a~calculation straightaway.
We need some preparations, in the course of which we will happen upon quaternion algebras.

\pagebreak[3]\txtskip

\begin{lemma}\label{lem:<v1,v2,v3>u=...}
Let\/ $u,\tinysp v_1,\tinysp v_2,\tinysp v_3\in\pure{A}$.
If\/ $\bilin{u}{v_i}\neq 0$ for some\/ $i=1,\,2,\,3$, then
\begin{equation}\label{eq:<v1,v2,v3>u=...}
\trilin{v_1}{v_2}{v_3}\negtinysp\narrt\cdot u
	\Eq \bilin{u}{v_1}(v_2\narrdt\cross v_3)
		\wide+ \bilin{u}{v_2}(v_3\narrdt\cross v_1)
		\wide+ \bilin{u}{v_3}(v_1\narrdt\cross v_2)~.
\end{equation}
\end{lemma}

\vspace{-1ex}\interskip

\begin{proof}
Write $d\defeq\trilin{v_1}{v_2}{v_3}$,
and denote by $w$ the pure element represented by the right hand side of~\eqref{eq:<v1,v2,v3>u=...}.
Easy calculations show that $\bilin{w}{v_i}=\bilin{u}{v_i}\negtinysp\narrt\cdot d$ for $i=1,\,2,\,3$,
and that $u\narrdt\cross w=0$.
For every $v\in\pure{A}$ we have
$0 = (u\narrdt\cross w)\narrdt\cross v = \bilin{u}{v}w - \bilin{w}{v}u$,
that~is, $\bilin{w}{v}u = \bilin{u}{v}w$.
Taking $v=v_1,\,v_2,\,v_3$,
we get $\bilin{u}{v_i}\negtinysp\narrt\cdot d\halftinysp u = \bilin{u}{v_i}\negtinysp\narrt\cdot w$
	for $i=1,\,2,\,3$,
and it follows that $d\halftinysp u = w$
because by assumption at least one $\bilin{u}{v_i}$ is non-zero.
\end{proof}

\thmskip

\begin{lemma}\label{lem:if-<v1,v2,v3>neq0...}
Suppose that\/ $v_1,\tinysp v_2,\tinysp v_3\in\pure{A}$ are such that\/ $\trilin{v_1}{v_2}{v_3}\neq 0$.
Then\/ $\tuple{v_1,v_2,v_3}$ is a~basis of\/~$\pure{A}$.
For~every non-zero\/ $u\narrt\in\pure{A}$
at least one of\/ $\bilin{u}{v_1}$, $\bilin{u}{v_2}$, $\bilin{u}{v_3}$~is\/~$\neq0$.
\end{lemma}

\interskip

\begin{proof}
We know that $v_1$, $v_2$, $v_3$ are linearly independent
because $d\defeq\trilin{v_1}{v_2}{v_3}\neq 0$.
Similarly $w_1\defeq v_2\narrdt\cross v_3$, $w_2\defeq v_3\narrdt\cross v_1$,
	$w_3\defeq v_1\narrdt\cross v_2$
are linearly independent because $\trilin{w_1}{w_2}{w_3}=d^{\tinysp2}$.
Clearly $\bilin{v_i}{w_j}=\delta_{ij}\tinysp d$ for $i,\tinysp j=1,\,2,\,3$.
Let\/ $V\negdtinysp$ be the subspace~of~$\pure{A}$ spanned by the $v_i$'s,
and let\/ $W\negdtinysp$ be the subspace of $\pure{A}$ spanned by the $w_i$'s.
Consider an arbitrary $u\in\pure{A}$.
If $\bilin{u}{v_1}\neq 0$, then $u\in W\negdtinysp$ by Lemma~\ref{lem:<v1,v2,v3>u=...}.
If $\bilin{u}{v_1}=0$, then $\bilin{u\narrdt+w_1}{\tinysp v_1}=d\neq 0$,
thus $u+w_1$ lies in $W\negdtinysp$, and so does~$u$.
It follows that $\pure{A}=W\negdtinysp$.
Since~$V\subseteq\pure{A}=W\negdtinysp$ and\/ $V\negdtinysp$ is three-dimensional,
we must have $V=W=\pure{A}$.

Regarding the second statement of the lemma, consider any non-zero $u\in\pure{A}$.
Since $\pure{A}=W\negdtinysp$,
we have $u=\xi_1w_1+\xi_2w_2+\xi_3w_3$ for some scalars $\xi_1$, $\xi_2$, $\xi_3$,
which are not all zero because $u$ is non-zero.
But $\bilin{u}{v_i} 
	= \xi_i\negtinysp\bilin{w_i}{v_i} = \xi_i\tinysp d$
for~$i=1,\,2,\,3$,
so at least one $\bilin{u}{v_i}$ is non-zero.
\end{proof}

\pagebreak[3]
\thmskip

An equivalent formulation of the second statement of Lemma~\ref{lem:if-<v1,v2,v3>neq0...}:
if $u,\tinysp v_1,\tinysp v_2,\tinysp v_3 \in \pure{A}$
and $u\neq 0$ and $\bilin{u}{v_i}=0$ for each $i=1,\,2,\,3$,
then $\trilin{v_1}{v_2}{v_3}=0$.

\thmskip

\begin{lemma}\label{lem:<v1,v2,v3>u=...always}
The identity~\eqref{eq:<v1,v2,v3>u=...}
	holds for all\/ $u,\tinysp v_1,\tinysp v_2,\tinysp v_3\in\pure{A}$.
\end{lemma}

\interskip

\begin{proof}
If $u=0$, the identity is $0=0$.
If $u\neq 0$ and $\bilin{u}{v_i}=0$ for each $i=1,\,2,\,3$,
then $\trilin{v_1}{v_2}{v_3}=0$ and once more the identity is $0=0$.
If $\bilin{u}{v_i}\neq0$ for some~$i=1,\,2,\,3$,
then the identity holds by Lemma~\ref{lem:<v1,v2,v3>u=...}.
\end{proof}

\thmskip

%

\begin{proof*}[Proof of identity~(\ref{eq:grand-identity}).]
Write $w_1\defeq v_2\narrdt\cross v_3$, $w_2\defeq v_3\narrdt\cross v_1$,
	$w_3\defeq v_1\narrdt\cross v_2$.
We~expand the~determinant $D$ on the left hand side of~\eqref{eq:grand-identity} along the first row,
applying~\eqref{eq:<u1xv1,u2xv2>=...} to the cofactors of the row's entries as we go,
then after some rearranging we apply~\eqref{eq:<v1,v2,v3>u=...} with~$u_1$~in~place~of~$u$,
and we are almost there:
\begin{align*}
D &\Eq \bilin{u_1}{v_1}\narrdt\cdot\bilin{u_2\narrdt\cross u_3}{\tinysp w_1}
	\widet+ \bilin{u_1}{v_2}\narrdt\cdot\bilin{u_2\narrdt\cross u_3}{\tinysp w_2}
	\widet+ \bilin{u_1}{v_3}\narrdt\cdot\bilin{u_2\narrdt\cross u_3}{\tinysp w_3} \\[.5ex]
	&\Eq \bigbilin{u_2\narrdt\cross u_3}
		{\,\bilin{u_1}{v_1}\negtinysp w_1
			\narrt+ \bilin{u_1}{v_2}\negtinysp w_2
			\narrt+ \bilin{u_1}{v_3}\negtinysp w_3} \\[.5ex]
	&\Eq \bigbilin{u_2\narrdt\cross u_3}{\dtinysp\trilin{v_1}{v_2}{v_3}\negtinysp u_1} \\[.5ex]
	&\Eq \trilin{v_1}{v_2}{v_3}\narrdt\cdot\bilin{u_2\narrdt\cross u_3}{\tinysp u_1} \\[.5ex]
	&\Eq \trilin{u_1}{u_2}{u_3}\narrdt\cdot\trilin{v_1}{v_2}{v_3}~. \tag*{\qed}
\end{align*}
\end{proof*}

\vspace{-1.5ex}
\pagebreak[3]\bigskip

\section{Quadratic spaces of arbitrary dimension}
\label{sec:quad-spacs-of-any-dim}

\medskip

In the next three sections (after this one)
we shall carry out a~classification (of sorts) of quadratic extension algebras.
At the highest level
we~shall classify quadratic extension algebras according to the rank of their norm.
With every quadratic extension algebra~$A$
there are associated the quadratic space $\pair{A,\Qnorm_{\!A}}$
and its quadratic subspace $\pair{\pure{A},\qnorm_{\!A}}$.
Denoting~by~$\varepsilon$ the quadratic functional~$\xi\mapsto\xi^2$ on~$K$,
we have $N_{\!A} = \varepsilon\perp\qnorm_{\!A}$,
and hence $\rank\Qnorm_{\!A}=\rank\qnorm_{\!A}+1$.
We~shall say that the rank of $\Qnorm_{\!A}$ is the rank of~$A$ and write~it $\rank A\tinysp$,
and similarly, that the rank of $\qnorm_{\!A}$ is the rank of $\pure{A}$
and write it~$\dtinysp\rank\pure{A}\tinysp$.
We~shall find that a~quadratic extension algebra can have only one of the three ranks
$1$,~$2$,~or~$4$, and that the quadratic extension algebras of rank $4$
are precisely the quaternion algebras.

\txtskip

Surely you have noticed that
we are not restricting ourselves to just the finite-dimensional quadratic extension algebras:
they may be of any dimension (if they can), finite or infinite,
and the same is true of the associated quadratic spaces.
Consequently, we will require, for our classification task,
a~few simple facts about quadratic spaces of arbitrary dimension.

\txtskip

Let $V\negdtinysp$ be a~vector space over a~field~$K\negtinysp$, of finite or infinite dimension.

\txtskip

A~functional\/ $q\colon V\to K$ is said to be \notion{quadratic}
if there exists a~bilinear functional $b$ on $V\negdtinysp$
such that $q(x) = b(x,x)$ for every $x\in V\negdtinysp$.
Setting $B(x,y)\defeq\thalf\bigl(b(x,y)+b(y,x)\bigr)$,
we have $q(x)=B(x,x)$, where $B$ is a~symmetric bilinear functional;
$B$ is uniquely deter\-mined by $q$,
namely $B(x,y)=\thalf\bigl(q(x\narrdt+y)\narrt-q(x)\narrt-q(y)\bigr)$.%
\footnote{Now I must have written down this formula for the $N\negtinysp$-th time,
where $N$ is a~not-so-small natural number.
It~is no longer such an exhilarating experience as it was the first time round.}
Let\/~$E$ be a~basis of\/~$V\negdtinysp$.
Every $x\in V\negdtinysp$ has a~unique representation $x=\Sum_{e\in E}\negtinysp\xi_e\tinysp e$
with only finitely many coordinates $\xi_e\in K$ different from~$0$.
If $x=\Sum_e \xi_e\tinysp e$ and $y=\Sum_e\negtinysp \eta_e\tinysp e$,
then $B(x,y) = \Sum_{e,f} B(e,f)\tinysp\xi_e\tinysp \eta_f$,
where $B(e,f)=B(f\negtinysp,e)$ for all $e,\tinysp f\in E$.
Save for the symmetry, the coefficients $B(e,f)$ are completely arbitrary:
if $\beta_{e\negtinysp f}$ ($e,\tinysp f\in E$) are any scalars
such that $\beta_{e\negtinysp f}=\beta_{f\negtinysp e}$ for~all basis vectors $e$~and~$f$,
then, with $x$ and $y$ as above,
the formula $B(x,y)\defeq\Sum_{e,f}\negtinysp\beta_{e\negtinysp f}\tinysp\xi_e\eta_f$
defines a~symmetric bilinear functional~on~$V\negdtinysp$.


Given a~quadratic functional $q$
and the associated symmetric bilinear functional $B$,
both on $V\negdtinysp$,
we have a~\notion{quadratic space} $\tuple{V,q,B}$.
We define quadratic subspaces, orthogonality $x\perp y$,
orthogonal sums of subspaces $U\negdtinysp\perp W\negdtinysp$, orthogonal complements $U^\perp\negtinysp$,
and so on,
just as for the finite-dimensional quadratic spaces.
The radical of the~quadratic space $V\negdtinysp$ is also defined in the same way,
namely $\rad V\defeq V^\perp\negtinysp$,
and $V\negdtinysp$ is said to be regular if $\rad V=0$.
We define the rank of the quadratic space $V\negdtinysp$ to be the codimension of its radical:
$\rank V \defeq \codim\rad V \eq \dim(V/\negtinysp\rad V\tinysp)$.
The most significant difference~between
finite-dimensional and infinite-dimensional quadratic spaces is this:
for a~regular~$V\negdtinysp$,
the linear mapping $V\mapsto\dualsp{V} : x\mapsto B(x,\anon)$
is~an~isomorphism of vector spaces if $V\negdtinysp$ is finite-dimensional,
while it is only injective if $V\negdtinysp$ is infinite-dimensional.%
\footnote{For any vector space $V\negdtinysp$ over $K$
	the dual space of $V\negdtinysp$ is denoted by $\dualsp{V}$.
If the dimension $m$ of $V\negdtinysp$ is infinite,
then the dimension of $\dualsp{V}$ is $\card{K}^m\geq 2^m>m$,
therefore $V\negdtinysp$ is not isomorphic to $\dualsp{V}$.}

\txtskip

We will need a~few facts about quadratic space of an arbitrary dimension.%
\footnote{We shall actually use only
Proposition~\ref{prop:orth-complem-of-findim-reg-subspc} for $\dim{U}=1$
and Proposition~\ref{prop:existence-of-findim-reg-subspcs} with $n=2$.}

\thmskip

\begin{proposition}\label{prop:max-reg-subspcs}
Let\/ $\tuple{V,q,B}$ be a quadratic space over\/ $K$.
If\/ $U$ is a maximal regular subspace of\/ $V\negdtinysp$,
then\/ $V=U\negdtinysp\perp U^\perp$ and\/ $U^\perp=V^\perp$.
If\/ $U$ is a~complementary subspace of the radical\/ $V^\perp$ in the space\/ $V\negdtinysp$,
then\/ $U$ is a~maximal regular subspace of\/~$V\negdtinysp$.
\end{proposition}

\interskip

\begin{proof}
Let $U$ be a maximal regular subspace of $V\negdtinysp$.
First, $U\inters U^\perp=\rad U=0$ because $U$ is regular.
In order to show that $U+U^\perp=V\negdtinysp$, consider any $v\in V\negdtinysp$.
If $v\in U$, we are done, so assume that $v\notin U$.
Since by maximality of $U$ the quadratic subspace $U+K\negtinysp v$ is not regular,
its radical contains a~non-zero vector $\upr$,
which lies in $U^\perp$ so it certainly does not belong to $U$;
but then $U+K\negtinysp v=U+K\negtinysp\upr$,
and we have $v\in U+K\negtinysp\upr\subseteq U+U^\perp$.
From $U\subseteq V\negdtinysp$ we get $U^\perp\supseteq V^\perp$.
As for the reverse inclusion, first notice that every non-zero $\upr\in U^\perp$ is isotropic,
since otherwise $U\negdtinysp\perp K\negtinysp\upr$ would be regular, contrary to maximality of~$U$;
but then $U^\perp$ is totally isotropic,
thus every $\upr\in U^\perp$ is orthogonal
to all vectors in $U^\perp$ as well as to all vectors in~$U$,
and it follows that $U^\perp\subseteq V^\perp$.

Now let $U$ be any~complementary subspace of $V^\perp$ in $V\negdtinysp$;
it is clear that $V=U\negdtinysp\perp V^\perp$.
Let $\Upr$ be the radical of the quadratic subspace~$U$.
The subspace $\Upr$ of $V\negdtinysp$ is orthogonal to $U$ as well as to $V^\perp$,
thus is orthogonal to $V\negdtinysp$, hence is contained in $V^\perp$,
and it follows that $\Upr=0$, which proves that $U$ is regular.
Suppose that a~regular subspace $U_1$ of $V\negdtinysp$ contains $U$.
Then $U_1=U\negdtinysp\perp(U_1\narrt\inters V^\perp)$
(since the lattice of subspaces of $V\negdtinysp$ is modular),
where $U_1\inters V^\perp$ is contained in the radical of the quadratic subspace $U_1$,
and it follows that $U_1\inters V^\perp=0$ and hence $U_1=U$;
this proves maximality of~$U$.
\end{proof}

\thmskip

\begin{corollary}\label{cor:rank=dim-max-reg-subspc}
The dimension of every maximal regular subspace of a~quadratic space
is equal to the rank of the quadratic space.
\end{corollary}

\thmskip

Proposition~\ref{prop:max-reg-subspcs}
not only proves that maximal regular subspaces of a~quadratic space~$V\negdtinysp$ exist,
it also characterizes them as the complementary subspaces of the radical $V^\perp$.
The existence part relies on Zorn's lemma (lurking in the background).
Suppose we have a~subspace $U$ of a $K\negdtinysp$-space $V\negdtinysp$
and want to prove that $U$ has a~complementary subspace~in~$V\negdtinysp$.
First we notice that the set $\coll{W}$ of all subspaces $W\negdtinysp$ of $V\negdtinysp$ such that $U\inters W=0$,
ordered by inclusion, is strictly inductive
(this is because the union of any nonempty chain in $\coll{W}$ belongs to~$\coll{W}$),
so it has maximal elements;
then it is easy to see that every maximal member of $\coll{W}$ is a~complementary subspace of~$U$.
Similarly we prove, using Zorn's lemma, that every regular subspace of a~quadratic space $V\negdtinysp$
is contained in a~maximal regular subspace.
To this end we only need to prove that
the union of a~nonempty chain~$\coll{U}$ of regular subspaces of~$V\negdtinysp$
is a~regular subspace of~$V\negdtinysp$:
every $u\in\rad\Union\coll{U}$ lies in some $U\in\coll{U}$ and belongs to $\rad U$,
hence $u=0$ since $U$ is regular.

\thmskip

\begin{proposition}\label{prop:orth-complem-of-findim-reg-subspc}
If\/ $U$ is a~finite-dimensional regular subspace of a~quadratic space\/ $\tuple{V,q,B}$ over\/ $K$,
then\/ $V=U\negdtinysp\perp U^\perp$.
\end{proposition}

\interskip

\begin{proof}
First we have $U\inters U^\perp=\rad U=0$.
In order to prove that $U+U^\perp=V\negdtinysp$,~con\-sider~any~$v\in V\negdtinysp$.
The mapping $U\to K : u\mapsto B(v,u)$ is a~linear functional on~$U$;%
~since $U$ is finite-dimensional and regular,
there exists $\vpr\in U$ so that $B(v,u)=B(\vpr,u)$ for every $u\in U$.
But~then $\vprpr\defeq v-\vpr\in U^\perp$, thus $v=\vpr+\vprpr\in U+U^\perp$.
\end{proof}

\thmskip

\begin{proposition}\label{prop:existence-of-findim-reg-subspcs}
If\/ $\tuple{V,q,B}$ is a~quadratic space over\/ $K$,
then for every natural number\/ $n\leq\rank V\negdtinysp$
there exists an\/~$n$-dimensional regular subspace of\/~$V\negdtinysp$.
\end{proposition}

\interskip

\begin{proof}
The proof is by induction on~$n$.
When $n=0$, the assertion is evidently true.
Now suppose $n>0$.
By induction hypothesis there exists an~$(n\narrdt-1)$-dimensional regular subspace $U$ of~$V\negdtinysp$.
We have $V=U\negdtinysp\perp U^\perp$ by Proposition~\ref{prop:orth-complem-of-findim-reg-subspc}.
Since $\dim U<n\leq \rank V\negdtinysp$,
the regular subspace $U$ is not a~maximal regular subspace of~$V\negdtinysp$,
according to Corollary~\ref{cor:rank=dim-max-reg-subspc}.
It~follows that $U^\perp$ is not totally isotropic,
since otherwise every subspace $U_1$ of $V\negdtinysp$ that properly contains $U$
would have a~radical containing $U_1\inters U^\perp\neq 0$
{\large(}this because we have $U_1=U\negdtinysp\perp (U_1\narrt\inters U^\perp)${\large)}.
Thus we can choose an anisotropic $\upr\in U^\perp$
and construct the $n$-dimensional regular subspace $U\negdtinysp\perp K\negtinysp\upr$ of~$V\negdtinysp$.
\end{proof}

\pagebreak[3]\bigskip

\section{Quadratic extension algebras of rank~1}
\label{sec:rank-1-quadext-algs}

\medskip

In this section
we examine the structure of a~quadratic extension $K\negdtinysp$-algebra $A$ of rank~$1$.

\txtskip

Since $\rank\qnorm_{\!A} = 0$,
the quadratic space $\bigpair{\pure{A},\qnorm_{\!A}}$ is totally isotropic.
This means that $\bilin{u}{v}=0$ for all pure $u$ and~$v$,
so~the product of any two pure elements is~pure,
hence the subspace $\pure{A}$ with the induced multiplication
is a~zero$\tinysp$-square algebra
(which, however, is not a~subalgebra of~$A$,
because it does not contain the multiplicative identity $1_A$ of~$A$).
The algebra $A$ is obtained from the zero$\tinysp$-square algebra $\pure{A}$
by adjoining the multiplicative identity to~it.
The zero$\tinysp$-square algebra $\pure{A}$ can be chosen arbitrarily:
adjoining an identity element to a~zero$\tinysp$-square algebra
always yields a~quadratic extension algebra of rank~$1$.

\txtskip

The subspace $\pure{A}$ is clearly an ideal of~$A$,
so we have the following result:
\vspace{1ex}
\begin{quote}
\textit{A~quadratic extension $K\negdtinysp$-algebra $A$ of rank\/~$1$ is simple \iff~$A=K$.}
\end{quote}
\vspace{1ex}
This is the only result about quadratic extension algebras of rank~$1$
that we will need later on to establish various characterizations of quaternion algebras.
However, we are curious about the structure of quadratic extension algebras of rank~$1$,
or, which is the same thing,
we~want to know how to construct arbitrary zero$\tinysp$-square algebras.

\txtskip

So let $S$ be any~zero$\tinysp$-square $K\negdtinysp$-algebra.
Since $x\tinysp y+yx = (x\narrt+y)^2-x^2-y^2 = 0$ for~all $x,\tinysp y\narrt\in S$,
the multiplication on $S$ is anticommutative,
thus the product $xy=\thalf(xy\narrt-yx)$ is at the same time a~Lie bracket.
The algebra $S$ is nilpotent; to be precise, $SSS=0$.
Indeed, any $x,\tinysp y,\tinysp z\in S$ satisfy the Jacobi's identity
\begin{equation*}
x\tinysp yz + yzx + zx\tinysp y \Eq 0~,
\end{equation*}
and since $y(zx)+(zx)y=0$, it follows that $x\tinysp yz=0$.

Let us define the \notion{annihilator}~of~$S$
as $\dtinysp\annih{S} \defeq \set{x\in S\suchthat xS=0}\tinysp$;
because of the anticommutativity of multiplication
we have also $\dtinysp\annih{S}=\set{x\in S\suchthat Sx=0}$.
It is clear that $\tinysp\annih{S}$ is a~subspace of the $K\negdtinysp$-space~$S$.
We can rewrite $SSS=0$ as $SS\subseteq\annih{S}$.

Set $W\defeq\annih{S}$, and let $V\negdtinysp$ be any complementary subspace of\/~$W\negdtinysp$
in~the $K\negdtinysp$-space~$S$.
We have an antisymmetric bilinear function
$V\negdtinysp\times V\negdtinysp\to W \narrt: \pair{v,\vpr}\mapsto v\tinysp\vpr$,
which completely determines the multiplication on~$S=V\negdtinysp\oplus W\negdtinysp$,
because for all $v,\tinysp\vpr\in V\negdtinysp$ and all $w,\tinysp\wpr\in W\negdtinysp$
we have $(v\narrt+w)(\vpr\narrt+\wpr) = v\tinysp\vpr$.
Moreover, $\set{v\in V\negdtinysp\suchthat vV=0}=0\tinysp$;
in particular, $V\negdtinysp$~is not one$\tinysp$-dimensional,
because then $\set{v\in V\negdtinysp\suchthat vV=0}=V\neq 0$.

Conversely, suppose we have $K\negdtinysp$-spaces $V\negdtinysp$ and $W\negdtinysp$
and an~antisymmetric bilinear function
$\mu\colon\negtinysp V\negdtinysp\narrdt\times V\negdtinysp\narrt\to W\negdtinysp$.
We define a~bilinear multiplication on $S\defeq V\negdtinysp\oplus W\negdtinysp$ by
\begin{equation*}
(v\narrt\oplus w)(\vpr\narrt\oplus\wpr) \Defeq 0\oplus\mu(v,\vpr)~.
\end{equation*}
The~$K\negdtinysp$-algebra $S$ is associative because $(SS)S=S(SS)=0$,
and it is a~zero$\tinysp$-square algebra because $\mu(v,v)=0$ for all $v\in V\negdtinysp$.
Moreover, $\annih{S} = \set{v\narrt\in V\negdtinysp\suchthat \mu(v,V)=0}\oplus W\negdtinysp$,
thus if we have chosen $\mu$
so that $\set{v\narrt\in V\negdtinysp\suchthat \mu(v,V)=0}=0$,
then $\tinysp\annih{S}=0\oplus W\negdtinysp$.

We will not go into detailed classification of zero$\tinysp$-square algebras,
	up to an~isomorphism.
This seems to be an open problem,
so it might be quite interesting to tackle it, but it would be of no use for our purpose.

\txtskip

A~final note.
Let $A\neq K\negtinysp$ be a~quadratic extension algebra of rank~$1$.
Then $\pure{A}\neq 0$ and also $W\negdtinysp\defeq\annih{\pure{A}}\neq\nolinebreak 0\tinysp$;
indeed, were $W\negdtinysp=0$, we would have $\pure{A}\dtinysp\pure{A}\subseteq W\negdtinysp=0$,
therefore $W\negdtinysp=\annih{\pure{A}}=\pure{A}\neq\nolinebreak 0$,~a~contra\-diction.
Let us split $\pure{A}\tinysp$ as $V\negdtinysp\oplus W\negdtinysp$.
If $\pure{A}\tinysp$ is not a~zero$\tinysp$-product algebra, then $\dim{V}\negdtinysp\geq 2$.
It is easy to see that it is possible to have $\dim{V}\negdtinysp=2$ and $\dim{W}\negdtinysp=1$,
and that in such a~case the $K\negdtinysp$-algebra $A$ is isomorphic to $\tQalgLam{0}{0}{K}$.

\pagebreak[3]\bigskip

\section{Quadratic extension algebras of rank~2}
\label{sec:rank-2-quadext-algs}

\medskip

In this section $A$ is a~quadratic extension $K\negdtinysp$-algebra of rank~$2$.

\txtskip

Now $\rank\qnorm_{\!A}=1$,
so there exists $u\in\pure{A}$ such that $\alpha\defeq u^2$ is a~non-zero scalar,
and $\pure{A} = K\negtinysp u\perp V\negdtinysp$, where $V = u^\perp = \rad\pure{A}$.
Every $v\in V\negdtinysp$ anticommutes with every $w\in\pure{A}$~be\-cause~$v\perp w$.
Since $(\xi u\narrt+v)^2=\xi^2u^2 + \xi(uv\narrt+vu)+v^2=\alpha\dtinysp\xi^2$
for all $\xi\in K$ and all~$v\in V\negdtinysp$,
we have $\sqr\pure{A}=\alpha\narrdt\cdot\sqr{K}$,
thus $A$ determines the element $\alpha\narrdt\cdot\sqr K^\times$
of the group $K^\times\negdtinysp/\tinysp\sqr K^\times$.

Let $v,\tinysp\vpr\in V\negdtinysp$.  We claim that $v\tinysp\vpr=0$.
First, $(v\tinysp\vpr)^2=v\tinysp\vpr\narrdt\cdot v\tinysp\vpr=-\vpr v^2 \vpr = 0$.
Since $u$ anticommutes with $v$ as well as with $\vpr$, $u$ commutes with $v\tinysp\vpr$.
Since $v\perp\vpr$, the product $v\tinysp\vpr$ is pure,
hence $v\tinysp\vpr=\xi u + w$ for some $\xi\in K$ and some $w\in V\negdtinysp$,
and we have $0 = (v\tinysp\vpr)^2 = \alpha\dtinysp\xi^2$,
thus $\xi=0$ and $v\tinysp\vpr = w$.
Since $u$ both commutes and anticommutes with $w$,
we have $uw=w\tinysp u=-uw$ and hence $uw=0$,
which implies that $w=0$ because $u$ is invertible,
and hence $v\tinysp\vpr=w=0$, as claimed.
The subspace $V\negdtinysp$ is certainly closed under multiplication, which is all-zero on~$V\negdtinysp$.

Let $v\in V\negdtinysp$.  We claim that $uv\in V\negdtinysp$ (and hence $vu=-uv\in V\tinysp$).
This is an immediate consequence of $u\perp v\tinysp$:
the product $uv$ is pure, and $u\perp uv$ because $uv=u\narrdt\cross v$.

We see that $V\negdtinysp$ is an ideal of $A$, thus $A$ can be simple only if $V\negdtinysp=0$.
Moreover, if~$V\negdtinysp=\nolinebreak0$, then $A=K+K\negtinysp u$ is simple
\iff\ $\alpha$ is a~non-square,
and in such a~case
$A$ is isomorphic to the quadratic extension field $K\bigl(\sqrt{\alpha\tinysp}\tinysp\bigr)$
of the field~$K$.

\medskip

\begin{widequote}
\textit{A quadratic extension $K\negdtinysp$-algebra of rank $2$
	is simple \iff\ it is a~quadratic extension field of~$K$.}
\end{widequote}

\medskip

As was the case with quadratic extension algebras of rank~$1$,
the result above is all we will need to know about quadratic extension algebras of rank~$2$
for the purpose of characterizations of quaternion algebras,
and as was the case with quadratic extension algebras of rank~$1$,
we nevertheless want to know more about the structure of quadratic extension algebras of rank~$2$.
This time we will pursue the classification all the way to the isomorphism classes;
though there's no real need to go into so much detail,
we will do it simply because we \emph{can} do it, and quite painlessly at that.

\txtskip

The linear transformation $\varphi\colon V\negdtinysp\to V\negdtinysp : v\mapsto u\tinysp v$
satisfies the identity $\varphi^2=\alpha\narrdt\cdot\id_V\negdtinysp$
because $u\tinysp(u\tinysp v)=u^2v=\alpha\narrdt\cdot v$ for every $v\in V\negdtinysp$.

Conversely, assume that we have a~vector space $A$ over~$K$,
two elements $1$~and~$u$~of~$A$ and a~subspace $V\negdtinysp$ of $A$
such that $A=K1\oplus K\negtinysp u\oplus V\negdtinysp$ (this is an internal direct sum),
and that we have also a~non-zero scalar $\alpha$
and a~linear transformation $\varphi\colon V\negdtinysp\to V\negdtinysp$.
Then there is one and only one bilinear multiplication $\pair{x,y}\mapsto x\tinysp y$ on $A$
which has the following properties:
\begin{itemize}
\item[\dia\:] $1\tinysp x = x1 = x$ for every $x\in A\tinysp$;
\item[\dia\:] $u^2 = \alpha\dtinysp$;
\item[\dia\:] $u\tinysp v = \varphi\tinysp v$ and $v\tinysp u = -\tinysp\varphi\tinysp v$
	for every $v\in V\negtinysp$;
\item[\dia\:] $v\tinysp\vpr = 0$ for all $v,\tinysp\vpr\in V\negdtinysp$.
\end{itemize}
If, in addition, $\varphi$ satisfies the identity $\varphi^2=\alpha\narrdt\cdot\id_V$,
then the vector space $A$, equipped with the multiplication,
is~a~quadratic extension $K\negdtinysp$-algebra of rank~$2$
with $\pure{A}=K\negtinysp u\oplus V\negdtinysp$.
Since $(\xi u\narrt+v)^2\negdtinysp=\xi^2\alpha\in K$
for all $\xi\narrt\in K$ and $v\narrt\in V\negdtinysp$,
it remains to prove associativity of multiplication.
It suffices to verify that $(w_1w_2)\tinysp w_3=w_1(w_2w_3)$,
where each $w_i$ is either~$u$ or belongs to~$V\negdtinysp$.
If $w_1=w_2=w_3=u$, then $(u\tinysp u)\tinysp u=\alpha\tinysp u=u\tinysp(u\tinysp u)$.
If at least two of~$w_i$'s belong to $V\negdtinysp$, then $(w_1w_2)\tinysp w_3=0=w_1(w_2w_3)$.
For~the~remaining three cases, let~$v\in V\negtinysp$:
then $(u\tinysp u)\tinysp v = \alpha\tinysp v=\varphi^2 v = \varphi(\varphi\tinysp v)
	= u\tinysp(u\tinysp v)$,
and $(u\tinysp v)\tinysp u = -\tinysp\varphi(\varphi\tinysp v)
	= \varphi(\negtinysp-\tinysp\varphi\tinysp v) = u\tinysp(v\tinysp u)$,
and $(v\tinysp u)\tinysp u = -\tinysp\varphi(\negtinysp-\tinysp\varphi\tinysp v)
	= \varphi^2 v = \alpha\tinysp v = v\tinysp(u\tinysp u)$.

\txtskip

In the forthcoming classification of quadratic extension algebras we distinguish two cases,
depending on whether $\alpha$ is a~square or not.

\txtskip

\textsc{Case~1:}\, $\alpha=\delta^{\tinysp2}\in\sqr K^\times$.

Let us say that a~$K\negdtinysp$-algebra is  of \notion{type~1}
if it is a~quadratic extension $K\negdtinysp$-algebra of rank~$2$
that satisfies the condition of this case.
(Recall that $\alpha\narrdt\cdot\tinysp\sqr\negtinysp K^\times$
is determined by~a~quadratic extension algebra of rank~$2$.)\,
Every~$K\negdtinysp$-algebra isomorphic to a~$K\negdtinysp$-algebra of type~1 is itself of type~1.

Replacing $u$ by $u/\delta$, we can assume that $u^2=1$, and hence that $\varphi^2=\id_V$.
The linear transformations
$\pi_+\defeq\thalf(\id_V\narrt+\varphi)$ and $\pi_-\defeq\thalf(\id_V\narrt-\varphi)$ of $V\negdtinysp$
are complementary projectors onto the invariant subspaces
$V_+\defeq\pi_+V\negdtinysp$ and $V_-\defeq\pi_-V\negdtinysp$ of~$\varphi$,
so we have $V=V_+\oplus V_-$.
Since $\varphi\tinysp\pi_+=\pi_+$ and $\varphi\tinysp\pi_-=-\pi_-$,
the linear transformation $\varphi$ of~$V\negdtinysp$
restricts on $V_+$ to $\id_{V_+}$ and on $V_-$ to $-\id_{V_-}$,
that is, $\varphi\tinysp v=v$ for every $v\in V_+$ and $\varphi\tinysp v=-v$ for every $v\in V_-$.
Set $\cardnum{m}_+\defeq\dim V_+$ and $\cardnum{m}_-\defeq\dim V_-$.
We can always assume that $\cardnum{m}_+\geq \cardnum{m}_-$:
if it happens that $\cardnum{m}_+<\cardnum{m}_-$, we simply replace $u$ by $-u$.
The~pair of cardinal numbers $\pair{\cardnum{m}_+,\tinysp \cardnum{m}_-}$
	with $\cardnum{m}_+\geq \cardnum{m}_-$
is completely determined by~$A\tinysp$
(we shall say that the pair and $A$ are \notion{associated}),
because $A$ determines the (unordered) pair of subspaces $\set{V_+,V_-}\tinysp$:
if $\upr$ is any pure element of $A$ such that $(\upr)^2=1$,
then $\upr=\pm u+v$ for some $v\in V\negdtinysp$,
and we have $\bigset{(1\narrt+\upr)\rad\pure{A},\,(1\narrt-\upr)\rad\pure{A}\dtinysp} = \set{V_+,V_-}$.

It is clear that two $K\negdtinysp$-algebras of type~1 are isomorphic \iff\
the pairs of cardinal numbers associated with them are equal.

Moreover, every pair of cardinal numbers
	$\pair{\cardnum{m}_+,\tinysp \cardnum{m}_-}$, $\cardnum{m}_+\geq \cardnum{m}_-$,
is associated with some $K\negdtinysp$-algebras of type~1,
that is, it~determines an isomorphism class of $K\negdtinysp$-algebras.
Given $\pair{\cardnum{m}_+,\tinysp \cardnum{m}_-}$, we define a~multiplication
on $A\defeq K\oplus K\oplus K^{(\cardnum{m}_+)}\oplus K^{(\cardnum{m}_-)}$
by
\begin{multline*}
\bigtuple{\tau_1,\tinysp\xi_1,\tinysp v_1,\tinysp\vpr_1}
		\narrdt\cdot\bigtuple{\tau_2,\tinysp\xi_2,\tinysp v_2,\tinysp\vpr_2} \Eq \\[.5ex]
	\Bigtuple{\tau_1\tau_2\narrdt+\xi_1\xi_2,\:
		\tau_1\xi_2\narrdt+\tau_2\xi_1,\:
		\tau_1v_2\narrdt+\tau_2v_1\narrdt+\xi_1v_2\narrdt-\xi_2v_1,\:
		\tau_1\vpr_2\narrdt+\tau_2\vpr_1\narrdt-\xi_1\vpr_2+\xi_2\vpr_1}~.
\end{multline*}
Then $A$ is a~$K\negdtinysp$-algebra of type~1
	associated with the~pair~$\pair{\cardnum{m}_+,\tinysp \cardnum{m}_-}$.

Examples.
The $K\negdtinysp$-algebra
$\bigtxtmtx{\begin{smallmatrix}\!K & \!K\!\\[.2ex]\!0 & \!K\!\end{smallmatrix}}$
is of type~1, asociated with $\pair{1,0}$.
The $K\negdtinysp$-algebra $\tQalgLam{1}{0}{K}$
is of type~1, associated with $\pair{1,1}$.

\txtskip

\textsc{Case~2:}\, $\alpha\in K^\times\negtinysp\narrt\setdiff\sqr K^\times$.

We shall say that a~$K\negdtinysp$-algebra is  of \notion{type~2}
if it is a~quadratic extension $K\negdtinysp$-algebra of rank~$2$
satisfying the condition of this case.

Let\/ $\tinysp\coll{U}$
be the set of all subspaces $U\negdtinysp$ of the $K\negdtinysp$-space $V\negdtinysp$
such that $U\negdtinysp\inters\varphi\tinysp U=0$, partially ordered by inclusion.
It is easy to see
that the union of a~nonempty chain of subspaces belonging to $\tinysp\coll{U}$
is a~subspace belonging to~$\tinysp\coll{U}$,
thus $\tinysp\coll{U}$ possesses maximal elements. 

Let\/ $U\negdtinysp\in\tinysp\coll{U}$
be such that $W\defeq U\negtinysp+\varphi\tinysp U \neq V\negtinysp$;
then $U$ is not maximal in~$\tinysp\coll{U}$.
In order to prove this,
first note that $W\negdtinysp$ is an invariant subspace of\/~$\varphi$
{\large(}because $\varphi(\varphi\tinysp U)=\alpha\tinysp U =U${\large)},
then choose any $v\in V\narrdt\setdiff W\negdtinysp$;
we claim that $\varphi\tinysp v\notin W\negdtinysp+K\negtinysp v = W\negdtinysp\oplus K\negtinysp v$.
Suppose, to the contrary, that $\varphi\tinysp v=w+\xi v$ for some $w\in W\negdtinysp$ and some $\xi\in K$;
then $\alpha\dtinysp v = \varphi^2 v = \varphi\tinysp w + \xi\varphi\tinysp v
	= (\varphi\tinysp w+\xi w)\widet+\xi^2\tinysp v$
implies $\alpha\dtinysp v=\xi^2\tinysp v$, hence $\alpha=\xi^2$, a~contradiction.
We~have a~direct sum $U\oplus\varphi\tinysp U\oplus K\negtinysp v\oplus K\negtinysp\varphi\tinysp v$
of subspaces of\/~$V\negdtinysp$,
and it follows that\/ $U\negtinysp\narrt\oplus K\negtinysp v\widet\in\tinysp\coll{U}$.

Pick a~maximal $U\negdtinysp\in\tinysp\coll{U}\tinysp$;
then $V\negdtinysp=U\negdtinysp\oplus\varphi\tinysp U\negtinysp$.
Define the~$K\negdtinysp$-isomorphism $\psi$ from $U\negdtinysp\oplus U$
(an~external direct sum) to~$V\negdtinysp$
by $\psi(v,\vpr)=v+\varphi\tinysp\vpr$,
and define the~$K\negdtinysp$-automorphism $\varphipr$ of\/ $U\negdtinysp\oplus U\negdtinysp$
by $\varphipr(v,\vpr)=(\alpha\tinysp\vpr\negdtinysp,\tinysp v)$.
Then $\psi\tinysp\varphipr\negtinysp(v,\vpr) = \varphi\tinysp v+\alpha\tinysp\vpr
	= \varphi\dtinysp\psi(v,\vpr)$,
thus $\psi$ is an isomorphism
$\pair{U\negdtinysp\narrdt\oplus U,\dtinysp\varphipr}\to\pair{V,\tinysp\varphi}$
of $K\negdtinysp$-spaces endowed with linear transformations.

The dimension $\cardnum{m}$ of a~maximal subspace $U$ is determined by~$A\tinysp$:
if\/ $\dim A$ is finite, then\/ $\cardnum{m}=\thalf\dim A \widet- 1$,
and if\/ $\dim A$ is infinite, then\/ $\cardnum{m}=\dim A$.
An isomorphism class of $K\negdtinysp$-algebras of type~2
is associated with a~pair $\pair{\overbar{\alpha},\tinysp \cardnum{m}}$,
where $\overbar\alpha$ is a~nontrivial element of the group $K^\times\negdtinysp/\sqr K^\times$
and $\cardnum{m}$ is~a~cardinal number.
A~construction of a~$K\negdtinysp$\nobreakdash-algebra of type~2,
given any such pair $\pair{\overbar{\alpha},\tinysp \cardnum{m}}$,
where $\overbar{\alpha}=\alpha\narrdt\cdot\sqr K^\times$
for some non-square $\alpha\in K^\times$:
we~equip the vector $K\negdtinysp$-space
	$A\defeq K\oplus K\oplus K^{(\cardnum{m})}\oplus K^{(\cardnum{m})}$
with the multiplication
\begin{align*}
&\bigtuple{\tau_1,\tinysp\xi_1,\tinysp v_1,\tinysp\vpr_1}
		\narrdt\cdot\bigtuple{\tau_2,\tinysp\xi_2,\tinysp v_2,\tinysp\vpr_2} \Eq \\[.5ex]
&\quad
	\Bigl(\tau_1\tau_2\narrdt+\alpha\tinysp\xi_1\xi_2,\:
		\tau_1\xi_2\narrdt+\tau_2\xi_1,\:
		\tau_1v_2\narrdt+\tau_2v_1
			\narrdt+\alpha\tinysp\xi_1\vpr_2\narrdt-\alpha\tinysp\xi_2\vpr_1,\:
		\tau_1\vpr_2\narrdt+\tau_2\vpr_1\narrdt+\xi_1v_2-\xi_2v_1\Bigr)~;\!\!\!\!
\end{align*}
this $A$ is a~$K\negdtinysp$-algebra of type~2
belonging to the isomorphism class associated with~$\pair{\overbar{\alpha},\tinysp \cardnum{m}}$.

Example.  The $K\negdtinysp$-algebra $\tQalgLam{\alpha}{0}{K}$,
$\alpha$ a~non-square, is of type~2, associated with $\pair{\overbar{\alpha},1}$.

\pagebreak[3]\bigskip

\section{Some characterizations of quaternion algebras}
\label{sec:quat-algs}

\medskip

Let $A$ be a~quadratic extension $K\negdtinysp$-algebra of rank at least $3$.


Since $\rank\pure{A}\geq 2$,
there exists a~two$\tinysp$-dimensional regular subspace $U\negdtinysp$
of the quadratic space $\bigpair{\pure{A},\qnorm_{\!A}}$.
We can choose an orthogonal basis $\pair{e,f}$ of $U\negdtinysp$;
then $a\defeq e^2=-\bilin{e}{e}$ and $b\defeq f^2=-\bilin{f}{f}$ are non-zero scalars,
and $\bilin{e}{f}=0$ hence $e\negtinysp f=e\narrdt\cross\negtinysp f$ is pure.%
~We~have
\begin{equation*}
\trilin{e}{f\negdtinysp}{\tinysp e\negtinysp f}
	\Eq \bilin{e\narrdt\cross\negtinysp f\negtinysp}{\dtinysp e\narrdt\cross\negtinysp f}
	\Eq -(e\narrdt\cross\negtinysp f)^2
	\Eq e^2\negtinysp f^2 - \bilin{e}{f}^2
	\Eq a\halftinysp b
	\Neq 0~,
\end{equation*}
thus, according to~Lemma~\ref{lem:if-<v1,v2,v3>neq0...},
the pure elements $e$, $f\negtinysp$, $\tinysp e\negtinysp f$
form a~basis of the $K\negdtinysp$-space $\pure{A}$,
and we see that $A$ is isomorphic to the quaternion algebra $\tQalgLam{a}{b}{K}$.

\thmskip

Recall that, given any $a,\tinysp b\in \punct{K}\!$,
one way to define the quaternion algebra $\tQalgLam{a}{b}{K}$
is as the $K\negdtinysp$-algebra with the~$K\negdtinysp$-basis $\tuple{1,i,j,k}$,
where $i^2=a$, $j^2=b$, $ij=k$ and $ji=\nolinebreak-k$.
Assuming associativity of multiplication,
we derive the rest of the multiplication table, namely
$k^2=-ji\narrdt\cdot ij=-ji^2j=-i^2j^2=-a\halftinysp b$,
$ik=i\narrdt\cdot ij=i^2j=aj$, $ki=-ji\narrdt\cdot i=\nolinebreak-aj$,
and~similarly $jk=-b\tinysp i$, $kj=b\tinysp i\tinysp$;
direct calculations then show that the multiplication determined by the full multiplication table
is indeed associative.%
\footnote{We also obtain an associative algebra if one of the scalars $a$ and $b$ is~$0$,
or even if both of them are~$0$, though then $\tQalgLam{a}{b}{K}$ is not a~quaternion algebra%
\,---\,it is just a quaternion\emph{ish} algebra.}
A~quaternion $K\negdtinysp$\nobreakdash-algebra is defined as any $K\negdtinysp$-algebra
isomorphic to $\tQalgLam{a}{b}{K}$ for some $a,\tinysp b\in\punct{K}\!$.
In other words, $Q$ is a~quaternion $K\negdtinysp$-algebra
\iff\ $Q$ is a~$K\negdtinysp$-algebra, there exist non-zero scalars $a$ and $b$,
and there exists a~basis $\tuple{1,e,f\negtinysp,\tinysp g}$ of the $K\negdtinysp$-space~$Q$,
so that $e^2=a$, $f^2=b$, and $e\negtinysp f=-\negtinysp f\negtinysp e=g\tinysp$;
we shall call any such basis an~\notion{$\pair{a,b}$-basis~of\/~$Q$};%
%
\footnote{Mark that a~quaternion algebra may have many $\pair{a,b}$-bases
for particular scalars $a,\tinysp b\narrt\in\punct{K}\negdtinysp$,
and that it may have $\pair{a,b}$-bases for different pairs of scalars $\pair{a,b}$.
To give a~simple example for the latter,
if a~quaternion algebra has an $\pair{a,b}$-basis,
then it has an $\pair{\lambda^2\negtinysp a,\tinysp\mu^2b}$-basis
for any non-zero scalars~$\lambda$~and~$\mu$.}
also, we shall call $\tuple{e,f\negtinysp,g}$ an~$\tuple{a,b}$-basis of $\pure{Q}$.
A~quaternion algebra~$Q$ with an~$\pair{a,b}$-basis $\tuple{1,e,f\negtinysp,\tinysp g}$
is a~quadratic extension algebra with $\pure{Q}=K\negtinysp e+K\negdtinysp f+K\negtinysp g$,
because $(\xi e\narrt+\eta f\narrt+\zeta g)^2=a\tinysp\xi^2+b\tinysp\eta^2-a\halftinysp b\tinysp\zeta^2$
is a~scalar for any scalars $\xi$, $\eta$, and $\zeta$.

\pagebreak[3]
\thmskip

\begin{lemma}\label{lem:a-b-gens-make-quat-alg}
Suppose that a~non-zero $K\!$-algebra $Q$ is generated by elements $e$ and $f$
such that $e^2=a$ and $f^2=b$ are non-zero scalars and $e\negtinysp f=-\negtinysp f\negtinysp e$.
Then $Q$ is a~quaternion $K\!$-algebra,
and $\tuple{1,e,f\negtinysp,\tinysp e\negtinysp f}$ is its $\pair{a,b}$-basis.
\end{lemma}

\interskip

\begin{proof}
Since $Q$ is a~non-zero algebra, it contains the field od scalars $K=K1$.
It is clear that the $K\negdtinysp$-space $Q$
is spanned by $1$, $e$, $f\negtinysp$, and $g\defeq e\negtinysp f\negtinysp$,
so~it remains to prove that these four elements are linearly independent.
Given an invertible~$y\in Q$ we define a~linear transformation $\kappa_y$ of~$Q$
by $\kappa_y\tinysp x\defeq\thalf(x+y\halftinysp x\halftinysp y^{-1})$ for~$x\in Q$.\linebreak[3]
The~elements $e$ and $f$ are invertible
{\large(}$e^{-1}=a^{-1} e$, $f^{-1}=b^{-1}\negdtinysp f${\large)},
so we have the linear transformations $\kappa_e$ and $\kappa_{\negdtinysp f}$.
Since $e1e^{-1}=1$, $e\tinysp e\tinysp e^{-1}=e$,
$e\negtinysp f\negtinysp e^{-1}=-\negtinysp f\negtinysp e\tinysp e^{-1}=-\negtinysp f$,
$ege^{-1}=e(-\negtinysp f\negtinysp e)e^{-1}=-e\negtinysp f=-g$,
we have $\kappa_e 1=1$, $\kappa_e e=e$, and $\kappa_e f=\kappa_e g=0$,
and similarly $\kappa_{\negdtinysp f}1=1$, $\kappa_{\negdtinysp f}\negtinysp f=f$,
and $\kappa_{\negdtinysp f}e=\kappa_{\negdtinysp f}g=0$.
Let~$x=\alpha+\xi e+\eta f\negtinysp+\zeta g \in Q$
(where~$\alpha,\tinysp\xi,\tinysp\eta,\tinysp\zeta\in K$);
then $\kappa_e\tinysp x=\alpha+\xi e$ and $\kappa_{\negdtinysp f}\tinysp x=\alpha+\eta f\negtinysp$,
hence $\kappa_e\tinysp\kappa_{\negdtinysp f}\tinysp x=\alpha\tinysp$.%
~Setting\linebreak[3] $\tau\defeq\kappa_e\halftinysp\kappa_{\negdtinysp f}$\linebreak[3]
we have $\tau x=\alpha$, $\tau(e\tinysp x)=a\tinysp\xi$, $\tau(f\negtinysp x)=b\tinysp\eta$,
and $\tau(g\tinysp x)=-a\halftinysp b\tinysp\zeta$.
But~then $x=0$ implies $\alpha=\xi=\eta=\zeta=0$,
which means that $1$,~$e$,~$f\negtinysp$,~$\tinysp g$ are linearly independent.
\end{proof}

\thmskip

\begin{proposition}\label{prop:quat-alg-is-CSA}
Every~quaternion $K\negdtinysp$-algebra is a~simple algebra with center~$K$.
\end{proposition}

\interskip

\begin{proof}
Let $Q$ be a~quaternion $K\negdtinysp$-algebra.
There exist $a,\tinysp b\in\punct{K}$
such that $Q$ has an $\pair{a,b}$-basis $\tuple{1,e,f\negdtinysp,\tinysp g}$.
If the element $x=\alpha+\xi e+\eta f+\zeta g$ of~$Q$
belongs to~$\Center(Q)$,
then comparing the coefficients in $e\tinysp x=xe$ we find that $\eta=\zeta=0$,
and similarly $f\negtinysp x=x\negtinysp f$ gives us $\xi=\zeta=0$,
therefore $x=\alpha\in K$.

Let $J\neq Q$ be an~ideal of~$Q$,
and let $Q\to Q/J : x\mapsto\overbar{x}$ be the natural projection.
The~elements $\overbar{1}$, $\overbar{e}$, $\overbar{f}\negtinysp$, $\overbar{g}$
of the non-zero $K\negtinysp$-algebra $Q/J$
obey the same multiplication table as the elements $1$,~$e$,~$f\negtinysp$,~$\tinysp g$ of~$Q$,
thus, by Lemma~\ref{lem:a-b-gens-make-quat-alg},
$Q/J$ is a~quaternion $K\negdtinysp$-algebra
with an~$\pair{a,b}$-basis $\tuple{\overbar{1},\overbar{e},\overbar{f},\overbar{g}}$,
and it follows that $J=0$.
\end{proof}

\thmskip

We shall need the following fact, borrowed from the theory of central simple algebras
{\large(}in~order to~prove the implication (3)$\Implies$(2)
in Proposition~\ref{prop:characterizations-of-quat-algs} below{\large)}:
every element~$x$ of a~four-dimensional central simple $K\negdtinysp$-algebra
is a~zero of its reduced characteristic polynomial
$X^2\negtinysp-\tinysp t(x)X\negtinysp+\tinysp n(x)\in K[X]$,
where $t(x)\in K$ is the reduced trace~of~$x$ and $n(x)\in K$ is the reduced norm of~$x$.

\txtskip

Now it is easy to derive,
from the foregoing results, the two highlighted results in the preceding two sections,
and the overview of algebras of dimension at most~$3$ in Section~\ref{sec:examps-quadext-algs},
the following characterizations of quaternion algebras
(we~omit the evident proofs):

\thmskip

\begin{proposition}\label{prop:characterizations-of-quat-algs}
The following properties of a\/~$K\negdtinysp$-algebra\/ $A$ are equivalent:
\begin{itemize}
\item[\rm{(1)}\:] $A$ is a~quaternion algebra;
\item[\rm{(2)}\:] $A\neq K$ is a~simple quadratic extension algebra with center\/~$K\negtinysp$;
\item[\rm{(3)}\:] $A\neq K$ is a~simple algebra with center\/ $K$
	and\, $\dim_K\negdtinysp A\leq4\tinysp$;
\item[\rm{(4)}\:] $A=K\oplus V\negdtinysp$ is a~ground Clifford algebra
	and\, $\rank\tuple{V,\sqr_{V\negdtinysp}}\geq 2$;
\item[\rm{(5)}\:] $A=K\oplus V\negdtinysp$ is a~ground Clifford algebra of dimension at least\/ $3$
	and the~qua\-drat\-ic~space\/ $\tuple{V,\sqr_{V\negdtinysp}}$ is regular.
\end{itemize}
\end{proposition}

\vspace{-1.5ex}
\pagebreak[3]\bigskip

\section{(Anti)automorphisms of a quaternion algebra}
\label{sec:(anti)automorphs-of-quatalg}

\medskip

To begin with, we consider properties of an (anti)automorphism~$h$
of an arbitrary quadratic extension algebra~$A$.


If $\alpha\in K$,
then $h(\alpha) = h(\alpha1) = \alpha h(1) = \alpha$,
that is, $h$ fixes all scalars.
Clearly $h$ preserves squares (whether it is anti- or not):
that is, $h(x^2) = h(x)^2$ for every~$x\in A$.
For~every $u\in A\narrdt\setdiff K$ we have $h(u)\in A\narrdt\setdiff K$
because $h$ is a~bijection that fixes all scalars.
If~$u\in\pure{A}$, then $h(u)\in A\narrdt\setdiff K$ and $h(u)^2 = h(u^2) = u^2 \in K$,
thus $h(u)\in\pure{A\tinysp}$;
since also $h^{-1}(v)\in\pure{A}$ for every $v\in\pure{A}$,
we see that $h$ restricts to an automorphism $\pure{h}$
	of the $K\negdtinysp$\nobreakdash-space~$\pure{A}$.
The~$K\negdtinysp$\nobreakdash-linear transformation $\pure{h}$ 
preserves the scalar product $\bilin{\anon}{\anon}$ on~$\pure{A}\tinysp$:
for all $u,\tinysp v\in\pure{A}$,
$\bigbilin{\tinysp\pure{h}u}{\tinysp\pure{h}v} = -\half\bigl(h(u)h(v)\narrt+h(v)h(u)\bigr)
	= h\bigl(-\half(uv\narrdt+vu)\bigr) = h\bigl(\bilin{u}{v}\bigr)
	= \bilin{u}{v}$.
(This is also a~straightforward consequence of the fact that $h$ preserves all squares.)
In~other words,
$\pure{h}$ is a~bijective isometry of the quadratic space $\bigtuple{\pure{A},\qnorm_A}$,
\ie, it belongs to the orthogonal group $\Orth(\qnorm_A)$.

We denote by $\Aut(A)$ the group of all automorphisms of the $K\negdtinysp$-algebra $A$,
and by $\AugAut(A)$ the group consisting of all automorphisms and of all antiautomorphisms
of the $K\negdtinysp$-algebra $A$.
Let us allow writing the conjugation on $A$ as an in-line ``$\widet*$'',
so that ${*}\tinysp x$ is an alternative notation for $x^*$, for any $x\in A$.
Every $h\in\AugAut(A)$ commutes with the conjugation:
$h\tinysp{*}={*}\tinysp h$, that is, $h(x^*)=h(x)^*$ for every $x\in A$.
Indeed, if $x=\alpha+u$ {\large(}$\alpha\in K$, $u\in\pure{A}\dtinysp${\large)},
then $h(x) = \alpha + h(u)$ with $h(u)\in\pure{A}$,
and $h(x^*) = \alpha - h(u) = \bigl(\alpha\narrt+ h(u)\bigr)^{\negdtinysp*} = h(x)^*$.
We~write $h^*\defeq h\tinysp{*}={*}\tinysp h$ for every $h\in\AugAut(A)$,
and denote by $\Aut(A)^*$ the set of all antiautomorphisms of the $K\negdtinysp$-algebra~$A$.
If $h\in\Aut(A)$, then $h^*\in\Aut(A)^*$, and if $h\in\Aut(A)^*$, then $h^*\in\Aut(A)$.

When $A$ is commutative (there \emph{are} commutative quadratic extension algebras),
every antiautomorphism of $A$ is in fact an automorphism of~$A$,
and hence $\AugAut(A)=\Aut(A)$.
But, if $A$ is not commutative, then no antiautomorphism of~$A$ is an automorphism of~$A$,
(in~particular, the conjugation is not an automorphism),
the group $\Aut(A)$ is a~normal subgroup of index~$2$ of the group $\AugAut(A)$,
and $\Aut(A)^*$ is the only nontrivial coset of $\Aut(A)$ in $\AugAut(Q)$.

We can always recover $h\in\AugAut(A)$ from $\pure{h}$, as $h=\id_K\oplus\pure{h}$,
which means that the homomorphism of groups $\AugAut(A)\to\Orth(\qnorm_A) \wide: h\mapsto \pure{h}$
is injective.

If $h$ is an automorphism,
then $\pure{h}$ preserves the pure products and the mixed products:
$\pure{h}(u\narrdt\cross v)=\pure{h}u\cross\pure{h}v$ 
and $\bigtrilin{\tinysp\pure{h}u}{\pure{h}v}{\pure{h}w\negtinysp}=\negtinysp\trilin{u}{v}{w}$
for all $u,\tinysp v,\tinysp w\in\pure{A}$.
If $h$ is an anti-automorphism,
then $\pure{h}(u\narrdt\cross v)=-\tinysp\pure{h}u\cross\pure{h}v$ 
and $\bigtrilin{\tinysp\pure{h}u}{\pure{h}v}{\pure{h}w\negtinysp}=-\negtinysp\trilin{u}{v}{w}$
for $u,\tinysp v,\tinysp w\in\pure{A}$.

\pagebreak[3]
\txtskip

Let $Q$ be a~quaternion $K\negdtinysp$-algebra.
(Recall that $Q$ is not commutative.)

\txtskip

Now the mixed product is a~non-zero alternating trilinear form on~$\pure{Q}$.
If $\varphi$ is a~linear transformation of $\pure{Q}$
{\large(}that is, an endomorphism of the $K\negdtinysp$-space~$\pure{Q}\dtinysp${\large)},
then $\det\varphi$ is the unique scalar $\delta$
such that $\trilin{\varphi u}{\varphi v}{\varphi w}=\delta\narrt\cdot\trilin{u}{v}{w}$
for all $u,\tinysp v,\tinysp w\in\pure{Q}$.
The determinant of an isometry of the regular quadratic space $\bigtuple{\pure{Q},\qnorm_Q}$
can only be $1$~or~$-1$.
The~special orthogonal group $\SpecOrth(\qnorm_Q)=\Orth_+(\qnorm_Q)$
is defined as the subgroup~of~$\Orth(\qnorm_Q)$
consisting of all $\sigma\in\Orth(\qnorm_Q)$~with $\det\sigma=1$.
We write the set of all $\sigma\in\Orth(\qnorm_Q)$ with $\det\sigma=-1$ as~$\Orth_-(\qnorm_Q)$.
The subgroup $\SpecOrth(\qnorm_Q)$ of $\Orth(\qnorm_Q)$ is normal of index~$2$,
and $\Orth_-(\qnorm_Q)=-\tinysp\SpecOrth(\qnorm_Q)$ is
the only nontrivial coset of $\SpecOrth(\qnorm_Q)$ in $\Orth(\qnorm_Q)$.

If $h\in\Aut(Q)$,
then $\bigtrilin{\tinysp\pure{h}u}{\pure{h}v}{\pure{h}w\negtinysp}=\trilin{u}{v}{w}$
	for all $u,\tinysp v,\tinysp w\in\pure{Q}$,
thus $\det\pure{h}=1$,
and if $h\in\Aut(Q)^*$, then $\det\pure{h}=-1$.
The group homomorphism $\AugAut(Q)\to\Orth(\qnorm_Q) : h\mapsto\pure{h}$
restricts to a~group homomorphism $\Aut(Q)\to\SpecOrth(\qnorm_Q)$.
We already know that both these group homomorphisms are injective.
They are in fact group isomorphisms:

\thmskip

\begin{proposition}\label{prop:Aut(Q)->SO(Q)-is-surjective}
Let\/ $Q$ be a~quaternion\/ $K\negdtinysp$-algebra, and let\/ $\sigma\in\Orth(\qnorm_Q)$.
If\/~$\sigma\in\Orth_+(\qnorm_Q)$, then\/ $\id_K\negdtinysp\oplus\sigma\in\Aut(Q)$,
and if\/ $\sigma\in\Orth_-(\qnorm_Q)$, then\/ $\id_K\negdtinysp\oplus\sigma\in\Aut(Q)^*$.
\end{proposition}

\interskip

\begin{proof}
Let $\sigma\in\Orth_+(\qnorm_Q)$.
There exists an $\tuple{a,b}$-basis $\tuple{e,f\negtinysp,e\negtinysp f}$ of $\pure{Q}$
for some $a,\tinysp b\in\punct{K}\negtinysp$.
Since $\sigma$ preserves scalar products,
we have $(\sigma e)^2 = -\negtinysp\bilin{\sigma e}{\sigma e}
		= -\negtinysp\bilin{e}{e} = e^2 =\nolinebreak a$,
and~likewise $(\sigma\negtinysp f)^2 = f^2 = b\tinysp$;
also $\sigma e\perp \sigma\negtinysp f$,
because $\bilin{\sigma e}{\sigma\negtinysp f}=\bilin{e}{f}=0$.
It~follows that
$\bigtuple{\sigma e,\tinysp\sigma\negtinysp f\negtinysp,\tinysp (\sigma e)(\sigma\negtinysp f)}$
is an $\tuple{a,b}$-basis of~$\pure{Q}$.
The non-zero vectors $\sigma(e\negtinysp f)$~and~$(\sigma e)(\sigma\negtinysp f)$
both lie in the one-dimensional subspace
	$\set{\sigma e,\tinysp\sigma\negtinysp f}^\perp$ of $\pure{Q}$,
so are proportional.
Since $\bigtrilin{\sigma e}{\tinysp\sigma\negtinysp f\negtinysp}{\tinysp\sigma(e\negtinysp f)}
	= \det(\sigma)\trilin{e}{f\negtinysp}{e\negtinysp f\tinysp} = ab$
and
$\bigtrilin{\sigma e}{\tinysp\sigma\negtinysp f\negtinysp}{\tinysp(\sigma e)(\sigma\negtinysp f)}
	=\nolinebreak ab$,%
~the~proportional vec\-tors~$\sigma(e\negtinysp f)$ and $(\sigma e)(\sigma\negtinysp f)$
must be equal,
so
$\bigtuple{1,\tinysp\sigma e,\tinysp\sigma\negtinysp f\negtinysp,\tinysp\sigma(e\negtinysp f)}
	=\nolinebreak \bigtuple{1,\tinysp\sigma e,\tinysp\sigma\negtinysp f\negtinysp,
						\tinysp (\sigma e)(\sigma\negtinysp f)}$
is~an~$\tuple{a,b}$-basis of $Q$,
and we conclude that $\id_K\leftt\oplus\sigma$ is an automorphism of~$Q$.

Let $\sigma\narrt\in\Orth_-(\qnorm_Q)$.
Since $-\sigma$ preserves scalar products
and $\det(\narrt{-\sigma})=(-1)^3\det(\sigma)=\nolinebreak1$, we have $-\sigma\in\Orth_+(\qnorm_Q)$,
thus $h\defeq\id_K\oplus-\sigma\in\Aut(Q)$ by what we have proved above,
and it follows that $h^* = {*}\tinysp h = \id_K\leftt\oplus\sigma\in\Aut(Q)^*$.
\end{proof}

\thmskip

\begin{corollary}\label{prop:isomorphism-Aut(Q)->SO(Q)}
For every quaternion $K\negdtinysp$-algebra $Q$
the mapping $\AugAut(Q)\to\Orth(\qnorm_Q) :\linebreak[3] h\mapsto\pure{h}$
is~a~group~isomorphism that restricts to the isomorphism\/ $\Aut(Q)\to\SpecOrth(\qnorm_Q)$.
\end{corollary}

\thmskip

For every invertible $y\in Q$
we have the inner automorphism $\inner_y\colon Q\to Q : x\mapsto yxy^{-1}$
of the $K\negdtinysp$-algebra~$Q$.
The mapping $\inner\colon Q^\times\to\Aut(Q) : y\mapsto c_y$ is a~homomorphism of groups.
Below we shall show that the homomorphism $\inner$ is surjective and that $\ker\inner=K^\times$.

Let $w$ be an anisotropic vector of the quadratic space $\bigtuple{\pure{Q},\qnorm_Q}$.
Then $\pure{Q} = K\negtinysp w\perp w^\perp$,
where $w^\perp$ is a~two-dimensional hyperplane of~$\pure{Q}$.
The linear transformation $\refl_w$ of $\pure{Q}$,
defined by $\refl_w u \defeq -wuw^{-1}$ for $u\in\pure{Q}$,
is the~reflection along the vector $w$ across the hyperplane $w^\perp$;
indeed, clearly $\refl_w w=-w$,
and if $u\in w^\perp$, then $\refl_w u = (-wu)w^{-1} = uww^{-1} = u$.
The reflection $\refl_w$ is an isometry of $\bigtuple{\pure{Q},\qnorm_Q}$ of determinant~$-1$,
and $\id_K\oplus\refl_w$ is the antiautomorphism $*\tinysp\inner_w=\inner_w^*$ of $Q$.

\thmskip

\begin{lemma}\label{lem:every-automorph-of-Q-is-inner}
Every automorphism of a~quaternion\/ $K\negdtinysp$-algebra\/ $Q$ is inner.
\end{lemma}

\interskip

\begin{proof}
Let $h\in\Aut(Q)$; then $\pure{h}\in\Orth_+(\qnorm_Q)$.
Since $\id_Q$ is certainly an inner automorphism of~$Q$,
we can assume that $\pure{h}\neq\id_{\scriptpure{Q}}$.
According to the Cartan$\tinysp$-Dieudonn\'e theorem,
$\pure{h}$ is a~composite of two reflections
{\large(}since $\det\pure{h}=1$,
$\pure{h}$ cannot be a~single reflection or a~composite of three reflections{\large)},
thus there exist two anisotropic vectors~$v,\tinysp w\in\pure{Q}$
such~that~$\pure{h}=\refl_v\refl_w$.
But then $h = \id_K\oplus\refl_v\refl_w = \inner_v^*\tinysp\inner_w^*
	= \inner_v{*}\tinysp{*}\tinysp\inner_w = \inner_v\tinysp\inner_w = \inner_{vw}$.
\end{proof}

\thmskip

\begin{proposition}
If\/ $Q$ is a quaternion\/ $K\negdtinysp$-algebra,
then\/ $\Aut(Q) \widedt\isomorph Q^\times\!/K^\times$.
\end{proposition}

\interskip

\begin{proof}
The group homomorphism $\inner\colon Q^\times\to\Aut(Q)$ is surjective,
in view of Lemma~\ref{lem:every-automorph-of-Q-is-inner}.
It~is~clear that the kernel of $\inner$ contains $K^\times$.
Conversely, if $y\in\ker c$, then certainly $y\neq 0$,
and~$yx=xy$~for~every~$x\in Q$ hence $y\in\Center(Q)=K$.
We conclude that $\ker c = K^\times$,
and that the surjective homomorphism $\inner$ induces an~isomorphism $Q^\times\!/K^\times\to\Aut(Q)$.
\end{proof}

\bigskip

\section{Orthogonal bases of pure quaternions}
\label{sec:orth-pure-bases}

\medskip

Let $Q$ be a~quaternion $K\negdtinysp$-algebra
with an~$\pair{a,b}$-basis $\tuple{1,e,f\negtinysp,\tinysp e\negtinysp f}$.
Setting $g\defeq e\negtinysp f$,
the full multiplication table for the pure basis elements $e$, $f$, $g$ is
\begin{gather*}
e^2=a\,, \qquad f^2=b\,, \qquad g^2=-a\halftinysp b\,, \\[.5ex]
e\negtinysp f=-\negtinysp f\negtinysp e=g\,,
	\quad f\negtinysp g=-g\negtinysp f=-b\halftinysp e\,,
	\quad ge=-eg=-a\negtinysp f\,.
\end{gather*}
This multiplication table is somewhat lopsided;
it is well suited for some purposes since it involves only two scalars $a$ and $b$,
but it definitely lacks symmetry.

\txtskip

In~order to bring about the desired symmetry,
we consider an arbitrary orthogonal basis $\tuple{e,f\negtinysp,g}$ of $\pure{Q}$.
Writing $a\defeq e^2\in\punct{K}$ and $b\defeq f^2\in\punct{K}$,
we note that $\tuple{1,e,f\negtinysp,e\negtinysp f}$ is an $\pair{a,b}$-basis of~$Q$.
The orthogonal complement of $\set{e,f}$ in $\pure{Q}$ is a~one-dimensional subspace of $\pure{Q}$
which contains both $e\negtinysp f\neq 0$ and $g\neq 0$,
thus $e\negtinysp f=\gamma g$ for some $\gamma\in\punct{K}$;
it~is~clear that $\gamma$ may be any non-zero scalar.
Now $f\negtinysp g=-b\tinysp\gamma^{-1}e$ and $ge=-a\gamma^{-1}\negtinysp f$,
which suggests that we set $\alpha\defeq-b\tinysp\gamma^{-1}$ and $\beta\defeq-a\gamma^{-1}\negtinysp$,
and this indeed makes the multiplication table symmetric:
\begin{equation}\label{eq:alpha-beta-gamma-mult-table}
\begin{gathered}
e^2=-\beta\tinysp\gamma\,,
	\qquad f^2=-\alpha\tinysp\gamma\,,
	\qquad g^2=-\alpha\tinysp\beta\,, \\[.5ex]
e\negtinysp f=-\negtinysp f\negtinysp e=\gamma g\,,
	\quad f\negtinysp g=-g\negtinysp f=\alpha\tinysp e\,,
	\quad ge=-eg=\beta\negtinysp f\,.
\end{gathered}
\end{equation}
We shall call such a~basis $\tuple{1,e,f\negtinysp,\tinysp g}$
of a~quaternion algebra $Q$ an~\notion{$\tuple{\alpha,\beta,\gamma}$-basis of\/~$Q$}.
Note that the multiplication rules $e^2=-\beta\tinysp\gamma$, $f^2=-\alpha\tinysp\gamma$,
	$e\negtinysp f=-\negtinysp f\negtinysp e=\gamma g$
imply all other entries of the multiplication table
(which is not at all surprising, of course);
for~example,
from $\gamma^2g^2
	= e\negtinysp f\negtinysp\narrdt\cdot e\negtinysp f
	= -\negtinysp f\negtinysp e\narrdt\cdot e\negtinysp f
	= -e^2\negtinysp f^2
	= -\alpha\tinysp\beta\tinysp\gamma^2$
we get $g^2=-\alpha\tinysp\beta$,
while $\gamma\negtinysp f\negtinysp g
	= f\negtinysp\narrdt\cdot e\negtinysp f
	= -\negtinysp f\narrt{\narrdt\cdot}f\negtinysp e
	= -\negtinysp f^2 e
	= \alpha\tinysp\gamma\tinysp e$
gives us $f\negtinysp g=\alpha\tinysp e$.
Let us compute the mixed product $\trilin{e}{f\negtinysp}{\tinysp g}$:
\begin{equation*}
\trilin{e}{f\negtinysp}{\tinysp g}
	\widedt= \bilin{e\negtinysp f\negtinysp}{g}
	\widedt= \bilin{\gamma g}{g}
	\widedt= -\gamma g^2
	\widedt= \alpha\tinysp\beta\tinysp\gamma~.
\end{equation*}
The symmetry persists.

\pagebreak[3]
\txtskip

Given $\alpha,\tinysp\beta,\tinysp\gamma\in\punct{K}$,
we let the quaternion $K\negdtinysp$-algebra $\Qalg_K(\alpha,\beta,\gamma)$
be the $K\negdtinysp$\nobreakdash-algebra with the basis $\tuple{1,i,j,k}$
that satisfies the multiplication table~\eqref{eq:alpha-beta-gamma-mult-table}
{\large(}with $\tuple{i,j,k}$ in place of $\tuple{e,f\negtinysp,\tinysp g}$, of course{\large)}.
There is really no need
to prove associativity of the $K\negdtinysp$\nobreakdash-algebra $\Qalg_K(\alpha,\beta,\gamma)$,
because we can construct it as the $K\negdtinysp$-algebra $\tQalgLam{-\beta\gamma}{-\alpha\gamma}{K}$
in which we then choose the $\tuple{\alpha,\beta,\gamma}$-basis $\tuple{1,i,j,\tinysp\gamma^{-1}ij}$.
However, a~reader having a~streak of mathematical masochism in him/$\tinysp$her
is going to enjoy the following verification (by~raw~brute~force)
of~the~associativity of the $K\negdtinysp$-algebra $Q=\Qalg_K(\alpha,\beta,\gamma)$.
Choosing~$Q$ to be $K^4\negtinysp$, a~$K\negdtinysp$-space with the standard basis
$1_Q=\tuple{1,0,0,0}$, $i=\tuple{0,1,0,0}$, $j=\tuple{0,0,1,0}$, $k=\tuple{0,0,0,1}$,
we define the (bilinear) multiplication on~$Q$ by
\begin{align*}
\tuple{\tau_1,\xi_1,\eta_1,\zeta_1}\tuple{\tau_2,\xi_2,\eta_2,\zeta_2} \Defeq
	\bigl(& \tau_1 \tau_2 - \beta\tinysp\gamma\tinysp \xi_1 \xi_2
			- \alpha\tinysp\gamma\tinysp \eta_1 \eta_2
			- \alpha\tinysp\beta\tinysp \zeta_1 \zeta_2, \\
	      & \tau_1 \xi_2 + \tau_2 \xi_1
			+ \alpha\tinysp \eta_1 \zeta_2 - \alpha\zeta_1 \eta_2, \\
	      & \tau_1 \eta_2 + \tau_2 \eta_1
			+ \beta\tinysp \zeta_1 \xi_2 - \beta\tinysp \xi_1 \zeta_2, \\
	      & \tau_1 \zeta_2 + \tau_2 \zeta_1
			+ \gamma\tinysp \xi_1 \eta_2 - \gamma\tinysp \eta_1 \xi_2\bigr)~.
\end{align*}
Now we consider three arbitrary elements
$x_i\defeq\tuple{\tau_i,\xi_i,\eta_i,\zeta_i}$ ($i=1,\,2,\,3$) of~$Q$,
and compare $(x_1x_2)\tinysp x_3$ with $x_1(x_2\tinysp x_3)$;
we obtain a~quartet of identities in the polynomial ring
$\ZZ[\alpha,\beta,\gamma,
	\tau_1,\xi_1,\eta_1,\zeta_1,\tau_2,\xi_2,\eta_2,\zeta_2,\tau_3,\xi_3,\eta_3,\zeta_3]$,
where for the~few moments it takes to make the comparison we regard
$\alpha$, $\beta$, $\gamma$, $\tau_1$, $\xi_1$, \ldots, $\eta_3$, $\zeta_3$
as distinct formal variables.%
\footnote{We do not carry out the comparison ourselves%
\,---\,manual calculations are too prone to mistakes\,---\,%
instead we hand the task over to a~computer program,
such as \Mathematica,
which is better at handling symbolic expressions than we are.}

\txtskip

This is how the constructions $\tQalgLam{a}{b}{K}$ and $\Qalg_K(\alpha,\beta,\gamma)$
of quaternion algebras are related:
\begin{equation*}
\Qalg_K(\alpha,\beta,\gamma) \Isomorph \QalgLam{-\beta\tinysp\gamma}{-\alpha\tinysp\gamma}{K}
\end{equation*}
for all $\alpha,\tinysp\beta,\tinysp\gamma\in\punct{K}\negdtinysp$,
and, in the other direction,
\begin{equation*}
\QalgLam{a}{b}{K} \Isomorph \Qalg_K(-b,-a,1)
\end{equation*}
for all $a,\tinysp b\in\punct{K}$.

\pagebreak[3]\bigskip

\section{Arbitrary bases of pure quaternions}
\label{sec:eny-pure-bases}

\medskip

Let $Q$ be a~quaternion $K\negdtinysp$-algebra.
This time we write the pure product of the elements $u$, $v$ of $\pure{Q}$ as $u\crossQ v$.
Let $\tuple{v_1,v_2,v_3}$ be an~arbi\-trary basis~of~$\pure{Q}$,
and set $w_1\defeq v_2\narrdt{\crossQ} v_3$, $w_2\defeq v_3\narrdt{\crossQ} v_1$,
	$w_3\defeq v_1\narrdt\crossQ v_2\tinysp$;
we know that $\tuple{w_1,w_2,w_3}$ is a~basis of $\pure{Q}$.

Let us regard the elements of the $K\negtinysp$-space $K^3$ as columns;
that is, $\tuple{x_1,x_2,x_3}\in K^3$~is
just an in-line notation for the column $\txtmtx{\widet{x_1~x_2~x_3}}\transp$.
For $x,\tinysp y\in K^3$ we let $x\narrt\cross y$ denote the usual cross product of triples.
We define isomorphisms of $K\negtinysp$-spaces $\varphi,\tinysp\psi\colon K^3\to\pure{Q}$~by
\begin{equation*}
\varphi\tinysp x \Defeq x_1v_1+x_2v_2+x_3v_3~, \qquad \psi\tinysp x \Defeq x_1w_1+x_2w_2+x_3w_3~,
\end{equation*}
for every $x=\tuple{x_1,x_2,x_3}\in K^3$.
Expanding vectors $w_1$, $w_2$, $w_3$ in the basis $\tuple{v_1,v_2,v_3}$
as $w_j=\Sum_{i=1}^3\negtinysp \alpha_{ij}\tinysp v_i$ ($j=1,\,2,\,3$)
and denoting the matrix $\txtmtx{\widet{\alpha_{ij}}}_{3\times 3}$ by $M$,
we have
\begin{equation*}
\psi\tinysp x \Eq \varphi Mx \qquad\quad \text{for every $x\in K^3$}\,.
\end{equation*}
Let $\Mpr\defeq\txtmtx{\narrt{\bilin{v_i}{v_j}}}_{3\times 3}$.
Then
$\bilin{\varphi\tinysp x}{\varphi\tinysp y} = x\transp\negtinysp\Mpr y$ for all $x,\tinysp y\in K^3$.

Let $d\defeq\trilin{v_1}{v_2}{v_2}$; then $\det\Mpr = d^{\tinysp2}$, by the `grand identity'.
The nine identities
\begin{equation*}
d\tinysp\delta_{ij} \Eq \bilin{v_i}{w_j} \Eq \Sum_{k=1}^3\narrt{\bilin{v_i}{v_k}}\alpha_{kj}~,
	\qquad\quad i,\, j = 1,\,2,\,3~,
\end{equation*}
can be rewritten as a~single identity $\Mpr\negtinysp M=d\tinysp I$,
and we see that $M=d\narrdt\cdot(\leftt\Mpr)^{-1}$ is an invertible symmetric matrix.
From $\det(\Mpr\tinysp)\det(\leftt M)=\det(d\tinysp I) = d^{\tinysp3}$ we get $\det M=d$,
whence $\Mpr=dM^{-1}=\adjugM$
(where $\adjugM$ is the adjugate of the matrix~$M$), thus
\begin{equation*}
\bilin{\varphi\tinysp x}{\varphi\tinysp y} \Eq x\transp\negtinysp\adjugM y
	\qquad\quad \text{for all $x,\tinysp y\in K^3$}~;
\end{equation*}
knowing the matrix $M$, we of course know the adjugate matrix $\adjugM$,
and we can compute the scalar product of the pure quaternions $\varphi\tinysp x$ and $\varphi\tinysp y$.
We can also compute the pure product, since 
\begin{equation*}
\varphi\tinysp x\crossQ\varphi\tinysp y
	\Eq \psi\tinysp(x\narrdt\cross y)
	\Eq \varphi\bigl(\narrt{M\negtinysp(\narrt{x\narrdt\cross y})}\bigr)
		\qquad\quad \text{for all $x,\tinysp y\in K^3$}~.
\end{equation*}
With the formulas for the scalar product and the pure product in our hands
we can now write out the formula for the product of general quaternions,~
\begin{equation*}
\bigl(s+\varphi\tinysp x\bigr)\bigl(t+\varphi\tinysp y\bigr)
	\Eq \bigl(s\tinysp t - x\transp\negtinysp\adjugM y\bigr)
		\widedt+ \varphi\bigl(s\tinysp y + tx + M\negtinysp(\narrt{x\narrdt\cross y})\bigr)~,
\end{equation*}
where $s,\tinysp t\in K$ and $x,\tinysp y\in K^3$.

\txtskip

For every nonsingular symmetric $3\narrdt\times 3$ matrix $M$ with entries in~$K$
there exists a~quaternion $K\negdtinysp$\nobreakdash-algebra~$Q$
with a~basis $\tuple{v_1,v_2,v_3}$ of $\pure{Q}$
such that $\psi=\varphi M$, where $\varphi,\tinysp\psi\colon K^3\to\pure{Q}$
are defined as above.
Given the matrix $M$,
we construct the quaternion $K\negdtinysp$-algebra $Q=\Qalg_K(\leftt M)=K^4=K\oplus K^3$
with the standard basis $\tuple{v_1,v_2,v_3}$ of $K^3=\pure{Q}$,
which~has the required property,
where the multiplication on $Q$ is defined by
\begin{equation*}
(s\oplus x)(t\oplus y)
	\Defeq \bigl(s\tinysp t - x\transp\negtinysp\adjugM y\bigr)
		\oplus \bigl(s\tinysp y+tx + M\negtinysp(x\narrdt\cross y)\bigr)~.
\end{equation*}
Associativity of multiplication is easily verified by raw brute force
(or via some clever shortcut).
Since $(0\oplus x)^2 = (-x\transp\negtinysp\adjugM x)\oplus 0$ for every $x\in K^3$,
we see that $Q$ is a~ground\linebreak[3] Clifford algebra.
And finally, $\rank\pure{Q}=3$ because $\adjugM$ is nonsingular,
thus $Q$ is indeed a~quaternion algebra.

The quaternion $K\negdtinysp$-algebra $Q$
with an $\tuple{\alpha,\beta,\gamma}$-basis, 
constructed in the preceding section,
is a~special case:
$\Qalg_K(\alpha,\beta,\gamma) = \Qalg_K\negtinysp\bigl(\diag(\alpha,\beta,\gamma)\bigr)$.

\bigskip

\section{Quaternion algebras over integral domains}
\label{sec:quatalgs-over-domains}

\medskip

This section offers a~glimpse
of the quaternion algebras over integral domains (of characteristic not~$2$).

\txtskip

Why under the sun would anyone want to observe a quaternion algebra $Q$
from a viewpoint of an arbitrary basis of $\pure{Q}\tinysp$?
There always exist orthogonal bases for the norm~$\qnorm_Q$,
and with respect to any such basis the matrix of the associated scalar product is diagonal,
which facilitates the exploration of the quaternion algebra.
However, we are able to diagonalize because the algebra $Q$ is over a~\emph{field}~$K$,
so that the quadratic functional~$\qnorm_Q$ is defined on the \emph{vector space} $\pure{Q}$.
What if we came across a~quaternion algebra defined over an~integral domain
(of characteristic not~$2$)?
Then it is no longer true that we can diagonalize at will
if a~quaternion algebra is presented relative to a~general basis;
when we do happen to find an orthogonal basis,
this appears to be nothing short of a~miracle
and is always an occasion (well, a~good excuse) for celebration.

\txtskip

The quaternion algebras over integral domains we are talking about are constructed as follows.
Let $R$ be an integral domain of characteristic different from~$2$,
and let $M$ be a~symmetric nonsingular $3\narrdt\times 3$ matrix with entries in~$R$.%
\footnote{Nonsingular $M$ means that $\det M\neq 0$
	(and not, perhaps, that $\det M$ is invertible in $R$).}
The quaternion algebra $Q=\Qalg_R(\leftt M)$
	is defined as the~free $R$-module $R\oplus R^{\tinysp3}$
{\large(}whose elements we write as~plain sums $r+u$,
	where $r\in R$ and $u\in R^{\tinysp3}${\large)}%
\footnote{That is,
we fearlessly identify $R\oplus 0$ with $R$ and $0\oplus R^{\tinysp3}$ with $R^{\tinysp3}$.}
with the $R$-bilinear multiplication defined by%
\footnote{This is the third\,---\,and the last\,---\,time
we write out the blasted formula for multiplication.}
\begin{equation}\label{eq:quatalg-multiplication}
\bigl(r\narrt+u\bigr)\bigl(s+v\bigr)
	\Defeq \bigl(rs-u\transp\negtinysp\adjugM v\bigr)
			+ \bigl(r\tinysp v + s\tinysp u + M\negtinysp(\narrt{u\narrdt\cross v})\bigr)~,
\end{equation}
for all $r\narrt+u,\tinysp s\narrt+v\in Q$.
We denote the vectors of the standard basis of $R^{\tinysp3}$ by $e_1$, $e_2$, $e_3$,
so that the general element of $Q$ is $s+u=s+u_1e_1+u_2e_2+u_3e_3$,
and denote by $\pure{Q}$ the direct summand $R^{\tinysp3}$ of~$Q$.

Let $K$ be the field of fractions of~$R$.
The $R$-algebra structure $Q=\Qalg_R(\leftt M)$ on $R\oplus R^{\tinysp3}$
uniquely extends to the $K$-algebra structure $Q_{\negdtinysp K}=\Qalg_K(\leftt M)$ on $K\oplus K^3$,%
\footnote{This is in effect the change of the base ring $R$ to the base field $K$:
	$Q_{\negdtinysp K}\isomorph K\narrt{\otimes_R} Q$.}
in which the multiplication is defined
by the same formula~\eqref{eq:quatalg-multiplication} as the multiplication in~$Q$,
only that now $r\narrt+u,\tinysp s\narrt+v\in Q_{\negdtinysp K}$.
The standard basis $\tuple{1,e_1,e_2,e_3}$ of the free $R$-module $Q$
is also a~basis of the $K\negdtinysp$-space~$Q_{\negdtinysp K}$.
Since $Q_{\negdtinysp K}$ is a~bona fide quaternion algebra over the field~$K$,
it makes perfect sense to call $Q$ a quaternion algebra over the integral domain~$R$.

We still have, on the three-dimensional free $R$-module $\pure{Q}=R^{\tinysp3}$,
the scalar product $\bilin{u}{v}=u\transp\negtinysp\adjugM v$ (with values in~$R\tinysp$),
the pure product $u\crossQ v=M(u\narrdt\cross v)$,
and the mixed product $\trilin{u}{v}{w}$.
Also there is the conjugation, which is an antiautomorphism of $Q$,
and there are the norm $\Qnorm\colon Q\to R$ and its restriction $\qnorm_Q:\pure{Q}\to R$.

\txtskip

Every automorphism $h$ of the $R$-algebra $Q$
extends to a~unique automorphism $h_K$ of the $K$-algebra $Q\tinysp$;
$h_K$ has the same matrix as $h$ with respect to the standard basis, and $h_K(Q)=Q$.
Conversely, if~$h$ is an automorphism of $Q_{\negdtinysp K}$
whose matrix with respect to the standard basis has all entries in~$R$
then $h$ restricts to the automorphism $h_R$ of $Q$;
indeed, since $\det h=1$, the matrix of $h$ is invertible in $\Malg_4(R)$,
thus the restriction $h_R$
is an automorphism of the $R$-module $Q$ and hence of the $R$-algebra~$Q$.

If $h$ is an automorphism of the $R$-algebra $Q$,
then its restriction $\pure{h}$ to $\pure{Q}$
is an automorphism of the quadratic $R$-module $\bigtuple{\pure{Q},\qnorm_Q}$ with $\det\pure{h}=1$,
that is, $\pure{h}\in\SpecOrth(\qnorm_Q)$.
Conversely, let $\sigma\in\SpecOrth(\qnorm_Q)$; then $\id_R\oplus\sigma\in\Aut(Q)$.
To prove this claim, we extend $\sigma$
to the automorphism $\sigma_K$ of the quadratic $K\negdtinysp$-space
	$\bigtuple{\pure{Q}_{\negdtinysp K},\rightt{\qnorm_{Q_{\negdtinysp K}}}}$;
we know that $h=\id_K\oplus\sigma_K$ is an automorphism of the $K\negtinysp$-algebra $Q_{\negdtinysp K}$,
and it is clear that the matrix of $h$ has all entries in $R$,
thus the restriction $h_R=\id_R\oplus\sigma$ is an automorphism of the $R$-algebra~$Q$.
It follows that the mapping $\Aut(Q)\to\SpecOrth(\qnorm_Q):h\mapsto\pure{h}$
is an isomorphism of groups, with the inverse $\sigma\mapsto\id_R\oplus\sigma$.

\txtskip

Every automorphism $h$ of the $K\negdtinysp$-algebra $Q_{\negdtinysp K}$ is inner:
there exists an invertible~$w$ in~$Q_{\negdtinysp K}$ such that $h(x)=wxw^{-1}$ for every~$x\in Q_{\negdtinysp K}$.
Rescaling $w$ by a~nonzero scalar in~$K$ yields the same automorphism~$h\tinysp$;
since $K$ is the field of fractions of $R$, we can assume that $w\in Q$.
Then $N(w)=ww^*\in\punct{R}$ and $h(x)=wxw^*\negdtinysp/\tinysp N(w)$,
and the matrix~$W\negdtinysp$ of the $K\negdtinysp$-linear transformation
$x\mapsto wxw^*$ of~$Q_{\negdtinysp K}$ has all entries in~$R$.
It follows that $h$~induces an automorphism of the $R$-algebra~$Q$
\iff\ all entries of the matrix~$W\negdtinysp$ are divisible by~$N(w)$.
But the matrix~$W\negdtinysp$ is of the form $W=\diag\bigl(N(w),\tinysp T\bigr)$,
where $T$ is the matrix
{\large(}relative to the standard basis $\tuple{e_1,e_2,e_3}$ of $Q_{\negdtinysp K}${\large)}
 of the $K\negdtinysp$-linear transformation $\pure{h}$ of~$\pure{Q}_{\negdtinysp K}$,
therefore $h$ induces an automorphism of $Q$
\iff\ every entry of the matrix $T$ is divisible by~$N(w)$.

Recall that $Q=\Qalg_R(\leftt M)$ and $Q_{\negdtinysp K}=\Qalg_K(\leftt M)$,
where $M\in\Malg_3(R)$ is a~nonsingular symmetric matrix which we write as
\begin{equation*}
M \Eq 	\begin{bmatrix}
	\,a_1 & b_3 & b_2\, \\
	\,b_3 & a_2 & b_1\, \\
	\,b_2 & b_1 & a_3\,
	\end{bmatrix}
	\,,
\end{equation*}
%
where the $a_i$ and the $b_i$ are elements of $R$ and
\begin{equation*}
H \Defeq \det M \Eq a_1a_2a_3 + 2\tinysp b_1b_2b_3-a_1b_1^2-a_2b_2^2-a_3b_3^2\Neq 0~.
\end{equation*}
We write the adjugate of the matrix $M$ as
\begin{equation*}
\adjugM \Eq	\begin{bmatrix}
		\,A_1 & B_3 & B_2\, \\
		\,B_3 & A_2 & B_1\, \\
		\,B_2 & B_1 & A_3\,
		\end{bmatrix}
		\,,
\end{equation*}
where $A_1 = a_1a_3-b_1^2$, $B_3=-a_3b_3+b_1b_2$, etc.
Let $w=w_0+w_1e_1+w_2e_2+w_3e_3\in\nolinebreak Q$ (the coordinates $w_k$ are in~$R$).
The norm of $w$ is
\begin{equation*}
N \Defeq N(w) \Eq
	w_0^2 + A_1w_1^2 + A_2w_2^2 + A_3w_3^2 + 2B_3w_1w_2 + 2B_2w_1w_3 + 2B_1\leftt{w_2}w_3~,
\end{equation*}
while the entries of the matrix $T$ (defined above) are
\begin{align*}
T_{11} &\Eq w_0^2 - 2\tinysp b_2w_0\tinysp w_2 + 2\tinysp b_3w_0\tinysp w_3
		+ A_1w_1^2 - A_2w_2^2 - A_3w_3^2
		- 2B_1\leftt{w_2}w_3~,\\[.5ex]
T_{12} &\Eq 2\tinysp b_2w_0\tinysp w_1 - 2\tinysp a_1\leftt{w_0}w_3
			+ 2B_3w_1^2 + 2A_2w_1w_2 + 2B_1w_1w_3~,
\end{align*}
and so on.

\thmskip

\begin{lemma}\label{lem:4(detM)wiwj=...}
Let\/ $M$, $H$, $\adjugM$, $w$, $N$, and\/ $T$ be as in the text.
Then:
\begin{align*}
4w_0^2 &\Eq T_{11} + T_{22} + T_{33} + N~, \\[.5ex]
4\rightt{H}w_0w_1 &\Eq B_2T_{12} - B_3T_{13} + B_1T_{22} - A_2T_{23} + A_3T_{32} - B_1T_{33}~, \\[.5ex]
4\rightt{H}w_0w_2 &\Eq -B_2T_{11} + A_1T_{13} - B_1T_{21} + B_3T_{23} - A_3T_{31} + B_2T_{33}~, \\[.5ex]
4\rightt{H}w_0w_3 &\Eq B_3T_{11} - A_1T_{12} + A_2T_{21} - B_3T_{22} + B_1T_{31} - B_2T_{32}~,
\end{align*}
\begin{align*}
4\rightt{H}w_1^2 &\Eq a_1T_{11} + 2\tinysp b_3T_{12} + 2\tinysp b_2T_{13} - a_1T_{22}
			- a_1T_{33} + a_1N~, \\[.5ex]
4\rightt{H}w_1w_2 &\Eq a_2T_{12} + b_1T_{13} + a_1T_{21} + b_2T_{23} - b_3T_{33} + b_3N~, \\[.5ex]
4\rightt{H}w_1w_3 &\Eq b_1T_{12} + a_3T_{13} - b_2T_{22} + a_1T_{31} + b_3T_{32} + b_2N~, \\[.5ex]
4\rightt{H}w_2^2 &\Eq -a_2T_{11} + 2\tinysp b_3T_{21} +a_2T_{22} + 2\tinysp b_1T_{23}
			- a_2T_{33} + a_2N~, \\[.5ex]
4\rightt{H}w_2w_3 &\Eq -b_1T_{11} + b_2T_{21} + a_3T_{23} + b_3T_{31} + a_2T_{32} + b_1N~, \\[.5ex]
4\rightt{H}w_3^2 &\Eq -a_3T_{11} - a_3T_{22} + 2\tinysp b_2T_{31} + 2\tinysp b_1T_{32}
			+ a_3T_{33} + a_3N~.
\end{align*}
\end{lemma}

\vspace{-1.5ex}
\thmskip

The identities in the lemma are verified by direct (if rather lengthy) calculations.%
\footnote{Don't do it by hand: just feed the (alleged) identities to \Mathematica.}
For the duration of verification we assume that the entries $a_i$ and $b_i$ of the matrix~$M$
and the coordinates $w_k$ of $w$ are formal variables;
then the identities live in the polynomial ring
$Z\defeq\ZZ[\tuple{a_i},\tuple{b_i},\tuple{w_k}]$.
This is proof by ``passing through the~generic~portal'':
the~lemma is true for every integral domain~$R$
\iff\ it is true for $R=Z$
{\large(}with~very special choices of $M\in\Malg_3(\leftt Z)$ and $w\in\Qalg_Z(\leftt M)${\large)}.
It suffices to verify only four identities%
\,---\,say the first, the second, the fifth, and the sixth identity\,---\,%
since the other six identities follow by symmetry (just rotate the indices $1$, $2$, $3$).

\txtskip

A \notion{GCD domain} is an integral domain in which every pair of elements has a~$\gcd$
(which is unique up to a~unit multiplier).
Let $R$ be a~GCD domain.
Then
\begin{equation*}
\gcd(ax_1,\ldots,ax_n) \Eq a\narrdt\cdot\gcd(x_1,\ldots,x_n)
\end{equation*}
(equality is up to a~unit factor) for all $a,\tinysp x_1,\tinysp\ldots,\tinysp x_n\in R$.%
\footnote{Much more is true.
If $K$ is a~field of fractions of a~GCD domain~$R$,
then the group $G\defeq K^\times\!/R^\times$,
ordered by divisibility over~$R$
(for $x,\tinysp y\in G$, $x\leq y$ iff $rx=y$ for some $r\in\punct{R}$),
is a~lattice-ordered~group.
See Chapter~2 in Stuart A.~Steinberg \cite{SteinbergLoRaM}.}
An element $x=\tuple{x_1,\ldots,x_n}$
of the free $R$-module $R^{\tinysp n}$ is said to be \notion{primitive}
if $x=r y$ with $r\in R$ and $y\in R^{\tinysp n}$ implies $r\in R^\times$;
$x$ is primitive iff $\gcd(x_1,\ldots,x_n)=1$.

\thmskip

\begin{proposition}\label{prop:Nw-divides-4detM}
Let\/ $R$ be a~GCD domain and\/ $K$ the field of fractions~of\/~$R$.
Let\/ $M$ be a symmetric nonsingular\/ $3\narrdt\times 3$ matrix with entries in\/ $R$,
and\/ $w$ a~primitive element of\/ $Q\defeq\Qalg_R(\leftt M)$ with\/ $N(w)\neq 0$.
If the automorphism\/ $x\mapsto wxw^{-1}\negdtinysp$ of the\/ $K\negdtinysp$-algebra\/~$Q_{\negdtinysp K}$
induces an automorphism of the\/ $R$-algebra\/~$Q$, then\/ $\Qnorm(w)$ divides\/~$4\det(\leftt M)$.
\end{proposition}

\interskip

\begin{proof}
The identities in Lemma~\ref{lem:4(detM)wiwj=...} tell us that $\Qnorm(w)$ divides
\begin{equation*}
\gcd\dtinysp\set{4\det(\leftt M)\tinysp w_iw_j\suchthat 0\leq i,\, j\leq 3}
	\Eq 4\det(\leftt M)\cdot\gcd\tinysp\set{w_iw_j\suchthat 0\leq i,\, j\leq 3}~,
\end{equation*}
where the gcd on the right hand side equals $\gcd(w_0,w_1,w_2,w_3)^2=1$.
\end{proof}

\pagebreak[3]
\thmskip

Let $R$ be an integral domain, $K$ its field of fractions, and $M\in\Malg_3(R)$ nonsingular symmetric.
The matrix $M$ determines a~``classic'' quadratic functional
	$q = q_M\colon R^{\tinysp3}\to R$,
where $q(x)=x\transp\negtinysp Mx$ for $x\in R^{\tinysp3}$.
We asociate with the quadratic $R$-module $\tuple{R^3\leftt,q}$
the quaternion $R$-algebra $Q\defeq\Qalg_R(\leftt M)$.
We do this because the~special orthogonal group of $q$
is isomorphic to the special orthogonal group of~$\qnorm_Q$.
In order to see this we identify $R$-linear transformations of $R^{\tinysp3}$
with their matrices relative to the standard basis;
then $\SpecOrth(q)$ is the set of all matrices $A\in\GenLin_3(R)$ with $\det A=1$
such that $A\transp\negtinysp M A=M$,
while $\SpecOrth(\qnorm_Q)$ is the set of all matrices $B\in\GenLin_3(R)$ with $\det B=1$
such that $B\transp\negtinysp\adjugM B=\adjugM$.
Inverting $A\transp\negtinysp M A=M$ {\large(}in $\Malg_3(K)${\large)}
and multiplying by $\det M$
we get $(A\invtransp)\transp\negtinysp\adjugM A\invtransp=\adjugM$;%
\footnote{Here $A\invtransp=(A^{-1})\transp=(A\transp)^{-1}$.
	If $A\in\GenLin_3(R)$, then $A\invtransp\in\GenLin_3(R)$;
the mapping $A\mapsto A\invtransp$ is an involution and an automorphism of $\GenLin_3(R)$.}
in the other direction we divide $B\transp\negtinysp\adjugM B=\adjugM$ by $\det M$
then invert, and obtain $(B\invtransp)\transp\negtinysp M B\invtransp=M$.
We therefore have the isomorphism of groups $\SpecOrth(q)\to\SpecOrth(\qnorm_Q):A\mapsto A\invtransp$
with the inverse $\SpecOrth(\qnorm_Q)\to\SpecOrth(q):B\mapsto B\invtransp$.

The most useful feature of $\SpecOrth(\qnorm_Q)$ is, of course, that it is isomorphic to $\Aut(Q)$.
Let~$B\in\SpecOrth(\qnorm_Q)$.
There exists $w\in\Qalg_R(\leftt M)$ with $\Qnorm(w)\neq 0$, such that $B=T_w/\Qnorm(w)$,
where $T_w$ is the matrix of the $R$-linear transformation $x\mapsto wxw^*$ of $\pure{Q}$
{\large(}and $\Qnorm(w)$ divides every entry of $T_w${\large)}.
Then $B\invtransp = T_{w^*}\transp/\Qnorm(w)\in\SpecOrth(q)$
{\large(}where $\Qnorm(w)$ divides every entry of the matrix $T_{w^*}\transp${\large)}.
Therefore, having a~general form of matrices in $\SpecOrth(\qnorm_Q)$,
we have also a~general form of matrices in $\SpecOrth(q)$.

We can explot the chain of isomorphisms
$\SpecOrth(q)\isomorph\SpecOrth(\qnorm_Q)\isomorph\Aut(Q)$
when we study the special orthogonal group~$\SpecOrth(q)$%
\,---\,its structure, its action on $R^{\tinysp3}$, etc.
Since there are only so many interesting results
that can be obtained for a~general integral domain~$R$,
the enquiry soon has to be specialized to, say, a~principal ideal domain~$R$,
or~even further to $R=\ZZ$ (the ultimately classical special case),
or~to $R=F[X]$ with~$F$ a~field of characteristic not~$2$ and $X$ a~formal variable.

And then a~serious fun begins\ldots\ 

But that is another story.

\bigskip

\section{Not an epilogue}
\label{sec:not-epilogue}

\medskip

We ran and ran after the stone that rolled and bounced down the slope.
It dislodged other stones, which in turn dislodged more stones\ldots\
so now we are chasing a~small avalanche of stones.
The bottom of the valley, obscured by roiling mists, is still very far below.
How and where and when will it all end?
Will it ever?


\end{document}

%% file: extended-preamble.tex
\input{general-preamble.tex}

\newcommand{\coll}{\mathcal}

\newcommand{\pair}{\arr}
\newcommand{\bigpair}{\bigarr}

\newcommand{\half}{\frac{1}{2}}
\newcommand{\thalf}{\tfrac{1}{2}}

\newcommand{\narr}[1]{\!#1\!}
\newcommand{\wide}[1]{\,#1\,}

\newcommand{\narrdt}[1]{\negdtinysp#1\negdtinysp}
\newcommand{\widedt}[1]{\dtinysp#1\dtinysp}
\newcommand{\narrt}[1]{\negtinysp#1\negtinysp}
\newcommand{\widet}[1]{\tinysp#1\tinysp}

\newcommand{\upr}{u^{\tinysp\prime}}
\newcommand{\vpr}{v^{\tinysp\prime}}
\newcommand{\wpr}{w^{\tinysp\prime}}
\newcommand{\xpr}{x^{\tinysp\prime}}

\newcommand{\vprpr}{v^{\tinysp\prime\prime}}

\newcommand{\Mpr}{M'}

\newcommand{\Upr}{U'}

\newcommand{\alphapr}{\alpha^{\tinysp\prime}}

\newcommand{\varphipr}{\varphi^{\tinysp\prime}}

%% file: general-preamble.tex
\usepackage{latexsym}
\usepackage{amssymb}
\usepackage{amsmath}
\usepackage{mathrsfs}  

\usepackage{amsthm}

\usepackage{graphicx}

\usepackage[margin=5mm,format=hang,justification=raggedright,font=small,labelfont=bf]{caption}

\usepackage{verbatim}

\usepackage{accents}

\usepackage{relsize}

\newcommand{\addtotextwidth}[1]{%
  \setlength{\hoffset}{.5\textwidth}
  \addtolength{\textwidth}{#1}
  \addtolength{\hoffset}{-.5\textwidth}}
\newcommand{\addtotextheight}[1]{%
  \setlength{\voffset}{.5\textheight}
  \addtolength{\textheight}{#1}
  \addtolength{\voffset}{-.5\textheight}}

\newlength{\itemsmargin}
\newlength{\itemslabelwidth}
\newlength{\itemslabelsep}
\newlength{\itemstopsep}
\newlength{\itemsitemsep}
  {\begin{list}{}{\setlength{\topsep}{\itemstopsep}
                  \setlength{\leftmargin}{\itemsmargin}
                  \setlength{\labelwidth}{\itemslabelwidth}
		  \setlength{\labelsep}{\itemslabelsep}
		  \setlength{\itemindent}{0pt}
                  \setlength{\listparindent}{.75\parindent}
                  \setlength{\itemsep}{\itemsitemsep}
                  \setlength{\parsep}{0pt}
                  \setlength{\parskip}{0pt}}}
  {\end{list}}

\renewenvironment{itemize}[1][2em]%
  {\begin{list}{}%
    {\setlength{\topsep}{1ex}
     \setlength{\labelwidth}{#1}
     \setlength{\leftmargin}{\labelwidth}
     \addtolength{\leftmargin}{1.5em}
     \setlength{\labelsep}{.618em}
     \setlength{\itemindent}{0pt}
     \setlength{\listparindent}{.75\parindent}
     \setlength{\itemsep}{1ex}
     \setlength{\parsep}{0pt}
     \setlength{\parskip}{0pt}}}
  {\end{list}}
  
\newcommand{\notion}[1]{{\bfseries #1}}

\renewcommand{\iff}{if~and only~if}

\newcommand{\ie}{i.e.}

\newcommand{\halftinysp}{\mspace{0.5mu}}
\newcommand{\tinysp}{\mspace{1mu}}

\newcommand{\dtinysp}{\mspace{2mu}}

\newcommand{\negtinysp}{\mspace{-1mu}}

\newcommand{\negdtinysp}{\mspace{-2mu}}

\newcommand{\MM}{\mathbb{M}}    

\newcommand{\ZZ}{\mathbb{Z}}    

\newcommand{\mathtextopfont}{\mathrm}
\newcommand{\mathtextop}[1]{\mathop{\smash{\mathtextopfont{#1}}}}

\newcommand{\union}{\cup}                       
\newcommand{\Union}{\bigcup}                    

\newcommand{\inters}{\cap}                      

\newcommand{\setdiff}{\mathbin{\smallsetminus}}

\renewcommand{\implies}{\mathrel{\Longrightarrow}}
\newcommand{\Implies}{\:\implies\:}

\newcommand{\isequiv}{\mathrel{\Longleftrightarrow}}
\newcommand{\Isequiv}{\:\isequiv\:}

\newcommand{\eq}{=}
\newcommand{\Eq}{\:=\:}
\newcommand{\Neq}{\:\neq\:}

\newcommand{\defeq}{\mathrel{\overset{\text{def}}{\tinysp=\!=\tinysp}}}
\newcommand{\Defeq}{\:\defeq\:}

\renewcommand{\leq}{\leqslant}

\renewcommand{\geq}{\geqslant}

\newcommand{\set}[1]{\{ #1 \}}                  
\newcommand{\bigset}[1]{\bigl\{ #1 \bigr\}}                  
\newcommand{\suchthat}{\mid}                    

\newcommand{\nothing}{}

\newcommand{\setparenarr}{\let\larr(\let\rarr)\let\arrmarg\nothing}
\newcommand{\setanglearr}{\let\larr\langle\let\rarr\rangle\let\arrmarg\tinysp}

\setparenarr

\newcommand{\arr}[1]{\larr\arrmarg #1 \arrmarg\rarr}         
\newcommand{\bigarr}[1]{\bigl\larr #1 \bigr\rarr}         
\newcommand{\Bigarr}[1]{\Bigl\larr #1 \Bigr\rarr}         
\newcommand{\tuple}{\arr}
\newcommand{\bigtuple}{\bigarr}
\newcommand{\Bigtuple}{\Bigarr}

\newcommand{\overlinepart}[2][1]%
  {{\setbox0=\hbox{$#2$}\accentset{\raisebox{0.0618ex}{\rule{#1\wd0}{0.09ex}}}{#2}}}
\newcommand{\overlinesymb}[1]{\overlinepart[0.8]{#1}}

\newcommand{\underlinepart}[2][1]%
  {{\setbox0=\hbox{$#2$}\underaccent{\rule{0pt}{.18ex}\rule{#1\wd0}{0.0618ex}}{#2}}}

\newcommand{\bigtruthord}[1]%
  {\text{\larger$\boldsymbol{[}$}\dtinysp{#1}\dtinysp\text{\larger$\boldsymbol{]}$}}

\newcommand{\abs}[1]{\mathopen{|}#1\mathclose{|}}           

\newcommand{\card}{\abs}                        

\newcommand{\anon}{\mathord{\halftinysp\rule[0.5ex]{0.5em}{0.5pt}\halftinysp}}

\newcommand{\Sum}{\sum}

\renewcommand{\dim}{\mathtextop{dim}\nolimits}		

\newcommand{\id}{\mathtextopfont{id}}

\newcommand{\farref}[1]{?.?}

\newcommand{\sgn}{\mathop{\mathrm{sgn}}\nolimits}

\newcommand{\twohyphs}%
	{\mathord{\dtinysp\raisebox{.15ex}{\text{-}}\dtinysp\raisebox{.15ex}{\text{-}}\dtinysp}}



%% file: preamble.tex
\newtheoremstyle{myplain}
  {0pt}
  {0pt}
  {\itshape}
  {}
  {\bfseries}
  {.}
  {0.75em}
  {}

\newtheoremstyle{mydefinition}
  {0pt}
  {0pt}
  {}
  {}
  {\bfseries}
  {.}
  {0.75em}
  {}

\makeatletter
\renewcommand*\env@matrix[1][*\c@MaxMatrixCols c]{%
  \hskip -\arraycolsep
  \let\@ifnextchar\new@ifnextchar
  \array{#1}}
\makeatother

\newenvironment{widequote}%
	{\begin{list}{}{%
		\setlength{\leftmargin}{\parindent}%
		\setlength{\rightmargin}{\parindent}}%
		\item[]\ignorespaces}%
	{\unskip\end{list}}

\theoremstyle{myplain}
\newtheorem{theorem}{Theorem} 
\newtheorem{proposition}[theorem]{Proposition}
\newtheorem{corollary}[theorem]{Corollary}
\newtheorem{lemma}[theorem]{Lemma}

\theoremstyle{mydefinition}

\newtheorem{remark}{Remark}

	{\begin{remark}}%
	{\hfill\qed\end{remark}}

\renewenvironment{proof}[1][Proof.]
	{\noindent{\bfseries #1}\hspace{.75em}\ignorespaces}
	{\hspace{\stretch{1}}\qed}
\newenvironment{proof*}[1][Proof.]
	{\noindent{\bfseries #1}\hspace{.75em}\ignorespaces}
	{\par}

\newcommand{\thmskip}{\bigskip\vspace{2ex}}
\newcommand{\txtskip}{\bigskip}
\newcommand{\interskip}{\bigskip}
\newcommand{\inskip}{\medskip}

\newcommand{\leftt}[1]{\negtinysp#1\tinysp}

\newcommand{\rightt}[1]{\tinysp#1\negtinysp}

\renewcommand{\emph}[1]{\textsl{#1\/}}

\renewcommand{\anon}{\mathord{\tinysp\rule[0.5ex]{0.5em}{0.5pt}\tinysp}}

\newcommand{\overbar}{\overlinesymb}

\newcommand{\dia}{\raisebox{.05ex}{\small{$\diamond$}}}
\renewcommand{\defeq}{:=}
\renewcommand{\Defeq}{\:\defeq\:}

\newcommand{\bilin}[2]{\left\langle#1,#2\right\rangle}
\newcommand{\bigbilin}[2]{\bigl\langle#1,#2\bigr\rangle}

\newcommand{\Mathematica}{\textsc{Mathematica}}
\newcommand{\diag}{\mathtextop{diag}\nolimits}
\newcommand{\adjug}{\widetilde}
\newcommand{\punct}[1]{#1^{\raisebox{0.25ex}{$\halftinysp\scriptscriptstyle\bullet$}}}

\newcommand{\Qnorm}{N\negtinysp}
\newcommand{\qnorm}{\nu}
\newcommand{\GenLin}{\mathrm{GL}}
\newcommand{\cross}{\times}
\newcommand{\Orth}{\mathrm{O}}
\newcommand{\SpecOrth}{\mathrm{SO}}

\newcommand{\pure}[1]%
	{\accentset{\raisebox{-0.25ex}[0pt][0pt]{$\smash{\scriptscriptstyle\rightharpoonup}$}}{#1}}
\newcommand{\Qalg}{\mathrm{Q}}
\newcommand{\transp}{^{\mathsf{T}}}
\newcommand{\invtransp}{^{-\mathsf{T}}}

\newcommand{\txtmtx}[1]{\boldsymbol{[}#1\boldsymbol{]}}
\newcommand{\bigtxtmtx}[1]%
	{\raisebox{-0.1ex}{\Large\boldmath$[$}#1\raisebox{-0.1ex}{\Large\boldmath$]$}}
\newcommand{\fieldchar}{\mathtextop{char}\nolimits}

\newcommand{\tQalgPete}[3]{\bigl(\negtinysp\narrt{\tfrac{\widet{#1,\dtinysp #2}}{#3}}\bigr)}
\newcommand{\QalgPete}[3]{\biggl(\negtinysp\narr{\frac{\widedt{#1,#2}}{#3}}\biggr)}
\newcommand{\tQalgLam}{\tQalgPete}
\newcommand{\QalgLam}{\QalgPete}

\newcommand{\Center}{\mathrm{Z}}

\newcommand{\Aut}{\mathrm{Aut}}
\newcommand{\AugAut}{\overbar{\mathrm{Aut}}}

\newcommand{\isomorph}{\cong}
\newcommand{\Isomorph}{\:\cong\:}
\newcommand{\Malg}{\MM}

\newcommand{\dualsp}[1]{#1^\wedge}

\newcommand{\scl}{\tau}
\newcommand{\vct}{\pi}
\newcommand{\trilin}[3]{\left\langle#1,#2,#3\right\rangle}
\newcommand{\bigtrilin}[3]{\bigl\langle#1,#2,#3\bigr\rangle}

\newcommand{\sqr}{\raisebox{-0.15ex}{$\Box$}\tinysp}

\newcommand{\rad}{\mathtextop{rad}\nolimits}

\newcommand{\rank}{\mathtextop{rank}\nolimits}
\newcommand{\cat}[1]{\mathbf{#1}}
\newcommand{\Cliff}{\mathit{Cl}}

\newcommand{\codim}{\mathtextop{codim}\nolimits}
\newcommand{\tr}{\mathtextop{tr}\nolimits}

\newcommand{\ExtAlg}{\Lambda}

\newcommand{\annih}{\mathtextop{ann}\nolimits}

\newcommand{\refl}{\varrho}
\newcommand{\inner}{c}

\newcommand{\adjugM}{\,\adjug{\!M}}

\newcommand{\crossQ}{\cross_Q}
\newcommand{\microrightharpoonup}{\includegraphics[hiresbb=true]{mp/micro-rightharpoonup-1.mps}}
\newcommand{\scriptpure}[1]%
	{\accentset{\raisebox{-0.25ex}[0pt][0pt]{\microrightharpoonup}}{#1}}
\newcommand{\cardnum}{\mathfrak}